\documentclass[11pt]{amsart}
\usepackage{amssymb}
\usepackage{amscd}
\usepackage{graphicx}
\usepackage[normalem]{ulem}

\usepackage[
hyperref,
]{degt}

\def\Dg:{\endgraf{\bf Dg:\enspace}\ignorespaces}
\def\2{\color{red}}
\def\3{\color{cyan}}
\def\4{\color{blue}}
\def\A:{\endgraf{\bf A:\enspace}\ignorespaces}

\def\bU{U_K}
\def\bV{\bar V}
\def\bT{\bar T}
\def\bX{X_K}

\def\lk{\operatorname{\ell\mathit{k}}}
\def\vlk{\operatorname{\,\overline{\!\lk}}}
\def\sg{\QOPNAME{sg}}

\def\dm{^*}

\def\Cs{\C^{\times\!}}

\def\CV{\Cal V}
\def\CVmax{\CV_{\max}}
\def\CA{\Cal A}
\def\CU{\Cal U}
\def\CCA{\CA^\circ}

\def\Clasp{\frak C}
\def\clasp{\frak c}
\def\open{^\circ}
\def\Fo{F\open}
\def\dL{\partial_L}
\def\dC{\partial_\Clasp}

\def\ZQ<#1,#2>{\{#1,#2\}}

\def\one{{\bar{\mathbf1}}}

\def\mone{{-\one}}
\let\mone=-

\def\CW{{\sl CW}}

\def\units{^{\times\!}}
\def\Cs{\C\units}
\def\int{\operatorname{int}}
\let\E=E
\def\D{\operatorname{D}}
\def\tGo#1{\widetilde\Go_{#1}}

\let\dif=R

\title{Slopes and concordance of links}

\author{Alex Degtyarev}

\address{%
Department of Mathematics\\
Bilkent University\\
06800 Ankara, Turkey}

\email{degt@fen.bilkent.edu.tr}

\author{Vincent Florens}

\address{%
Laboratoire de Math\'{e}matiques et leurs applications, UMR CNRS 5142\\
Universit\'{e} de Pau et des Pays de l'Adour\\
Avenue de l'Universit\'{e}\\
 BP 1155 64013 Pau Cedex, France
}

\email{vincent.florens@univ-pau.fr}

\author{Ana G.\ Lecuona}

\address{%
School of Mathematics and Statistics, University of Glasgow, Glasgow, UK.}

\email{ana.lecuona@glasgow.ac.uk}



\begin{document}

\begin{abstract}
The slope is an isotopy invariant of colored links with a distinguished
component, initially introduced by the authors to describe an extra
correction term in the computation of the signature of the splice. It
appeared to be closely related to several classical invariants, such as the
Conway potential function or the Kojima $\eta$-function (defined for
two-components links). In this paper, we prove that the slope is invariant
under colored concordance of links.
Besides,
we present a formula to compute the slope in terms of C-complexes
and generalized Seifert forms.
\end{abstract}

\maketitle

\section{Introduction}\label{S.intro}

The slope is an isotopy invariant defined for  so-called $(1,\mu)$-colored links
$K \cup L$ (with a distinguished component $K$ given color~$0$)
in rational homology spheres.
It is closely related to several classical invariants \cite{DFL1,DFL2,DFL3},
such as the
Conway potential
and Kojima-Yamasaki $\eta$-function (defined for two-components links
\cite{KY,Jin,Co}). To certain $\C\units$-valued characters
of the group $\pi_1(S\subset L)$,
\viz. those trivial on~$[K]$, see~\eqref{eq.CCA},
the slope associates a complex number (possibly infinite). The torus of
characters preserving the coloring is naturally identified with the complex
torus $(\C^\times)^\mu$,
and the slope is a function on (a Zariski open subset of) the variety
$\CA(K/L)\subset(\C^\times)^\mu$ of admissible characters. This function is
rational
away from a certain singular locus determined by the
Alexander module of $K \cup L$; however,
in general, the \emph{values} of the slope are not determined by the
Alexander module.

The aim of this paper is to show that the slope is invariant under colored
concordance of links, see \autoref{th.main}, and to present a method to
compute the slope in terms of the Seifert forms of the colored link $L$ with
an extra piece of data, see \autoref{th.C-complex}. In the case of
algebraically split links of two components, the invariance of the slope
under colored concordance was known for certain values, \viz. those where it:
coincides with the $\eta$-function \cite[Corollary~3.24]{DFL2}.
In this paper we show that, outside a certain subset of  $(\C^\times)^\mu$, the \emph{Knottennullstellen}
\cite{NP17,CNT}, concordant links have the same slope.
More generally, for algebraically split links with an arbitrary number of components, our result implies that a certain quotient of the Conway functions of $K \cup L$ and $L$ is invariant
 under colored concordance of $K \cup L$ (see \autoref{cor.nabla}), whereas the Conway functions themselves are \emph{not} concordance invariants (see \cite{Kau}).

One can compute the slope directly from the definition using the Fox
calculus; this is explained in \cite[Section~3.2]{DFL2}.
While
allowing for
easy computer assisted computations,
this approach
is not particularly useful when
dealing with families of examples.
In certain cases,
the slope can also be computed
as a quotient of the Conway polynomials (see \cite[Theorem~3.1]{DFL2}),
but this formula is inconclusive at the common roots of the
numerator and denominator (l'H\^{o}pital's rule does not work); in
particular, it leaves wide open the most interesting case where both
polynomials vanish identically.
In this
paper we
suggest
yet another method of computing the slope, using
C-complexes. These were introduced by Cooper \cite{Cooper} and extended, in
very recent years, by different groups to compute many link invariants
(\cite{Cimasoni1,CF,CFT,Merz} among others) and to study their properties
(\cite{DavisRoth,DavisTOtto,various} among others).

The computation of the slope using C-complexes is particularly powerful when
dealing with families of examples as in \cite[Example~3.28]{DFL2} and
\cite[Example~5.5]{DFL3}.
For the moment,
our formula only works in the special case of $K$ algebraically unlinked from
each monochrome component~$L_i$. For an
algebraically split two component link, the C-complex used in the
computation is merely
a Seifert surface.

The paper is organized as follows. In Section \ref{section.basics} we recall the construction and the basic
properties of the slope. Section \ref{section.conc} is devoted to the proof
of the
concordance invariance.
In Section \ref{section.comp} the
computation of the slope in terms of (generalized) Seifert forms is given,
and the main formula is proved in Section \ref{section.proof}.

\section{Slopes} \label{section.basics}

A \emph{$\mu$-colored link} is an oriented link~$L$ in $S^3$ equipped with a
surjective map $\pi_0(L)\onto\{1,\ldots,\mu\}$, called
\emph{coloring}.
The union of the components of~$L$ given the same color $i=1,\ldots,\mu$ is
denoted by $L_i$.
Each link has a canonical \emph{maximal coloring}, where each component
is given a separate color. In this special case, each $L_i$ is a knot.

We denote by $X:=S^3\sminus T_L$
the complement of a small open
tubular neighborhood of $L$.
The group $H_1(X)$ is free abelian, generated by the classes $m_C$
of the meridians of the components $C \subset L$.
By convention, $m_C$ is oriented so that $m_C\circ\ell_C=1$ in
$\partial T_C$, where $\ell_C$ is a longitude and the orientation on
$\partial T_C$ is that induced from~$X$.
The coloring induces an epimorphism
\[*
\varphi\: \pi_1(X) \onto H:=\bigoplus_{i=1}^\mu \Z t_i
\label{eq.H}
\]
sending $m_C$ to $t_i$ whenever $C \subset L_i$.
A multiplicative character $\Go\:\pi_1(X)\to\C\units$
is determined by its values on the meridians, and the torus of
characters preserving the coloring
(\ie, those that factor through~$\Gf$)
is naturally identified with the complex torus $(\C\units)^\mu$.
Through this identification, we set $\Go_{i}:=\Go(\varphi(t_{i}))$
and, with a certain abuse of the language,
speak about a character $\Go=(\Go_1,\ldots,\Go_\mu)$.
We define
\[*
\Go\1:=(\Go_1\1,\ldots,\Go_\mu\1),\qquad
\bar\Go:=(\bar\Go_1,\ldots,\bar\Go_\mu),\qquad
\Go^*:=\bar\Go\1.
\]
A character~$\Go$ is called \emph{unitary} if
$\Go^*=\Go$, \ie,
$\ls|\Go_i|=1$ for all $i=1,\dots,\mu$.
Unitary characters constitute a torus $(S^1)^\mu\subset(\C\units)^\mu$.

Given a topological space~$X$ and a multiplicative character
$\Go\:\pi_1(X)\to\C\units$, we denote by $H_*(X;\C(\Go))$ the homology
of~$X$ with
coefficient in the local system $\C(\Go)$ twisted by~$\Go$. See \cite[Section 2]{DFL2} for more details.

In this paper, we consider mainly colored links with a distinguished component. They are
 \emph{$(1,\mu)$-colored links}, defined as $(1+\mu)$-colored links of the form
\[*
K\cup L=K\cup L_1\cup\ldots\cup L_\mu,
\]
where the \emph{knot}~$K$ is the only component given the distinguished
color~$0$. The \emph{linking vector} of a $(1,\mu)$ colored link is
$\vlk(K,L):=(\Gl_1,\ldots,\Gl_\mu)\in\Z^\mu$, where $\Gl_i:=\lk(K,L_i)$.

\definition
A character $\Go\: \pi_1( X)\longrightarrow \C\units$
on a $(1,\mu)$-colored link $K\cup L$ is
called \emph{admissible} if $\Go([K])=1$;
it is called \emph{non-vanishing} if
$\Go_{i}\neq 1$ for all $i=1,\dots,\mu$.
\enddefinition

The variety of admissible
characters is denoted $\CA(K/L)$,
and $\CCA(K/L)\subset\CA(K/L)$ is the subvariety of admissible
non-vanishing characters.
Letting $\Gl:=\vlk(K,L)$, we have
\[
\CA(K/L) = \bigl\{ \Go \in (\C\units)^\mu\bigm|\Go^\Gl=1\bigr\},\quad
\CCA(K/L)=\CA(K/L)\cap(\C\units \sminus 1)^\mu,
\label{eq.CCA}
\]
where $\Go^\lambda:=\prod\Go_i^{\lambda_i}$.
In particular, if $\Gl=0$,
then $\CCA(K/L)=(\C\units\sminus1)^\mu$.

Let
$\bX = S^3\sminus T_{K \cup L}$ be the complement of an open tubular neighborhood of $K \cup L$.
We abbreviate $m:=m_K$ and $\ell:=\ell_K$, where $\ell_K$ is the
\emph{preferred} (\viz. unlinked with~$K$) longitude, also called
\emph{Seifert longitude}.

\remark \label{rem.charac}
Any character
$\Go\in(\C\units)^\mu$ extends to
a natural character
$ \pi_1(\bX) \rightarrow \C\units$ sending
$m$ to $1$; for short, this
extension
is also denoted $\Go$.
In this language,
the original character~$\Go$ is admissible if and only if $\Go(\ell)=1$.
\endremark

We denote by
$\partial_K\bX=\partial T_K$ the intersection of
$\partial\bX$ with the closure of $T_K$ and
consider the inclusion
\[*
i\: \partial_K \bX {}\into\partial\bX \hookrightarrow \bX.
\]

If
$\Go \in \CCA(K/L)$,
the induced homomorphism
\[
i_*\: H_1(\partial_K \bX; \C(\Go))
{}\overset\simeq\longto H_1(\partial\bX;\C(\Go))
\longto H_1(\bX;\C(\Go))
\label{eq.boundary}
\]
can be regarded as that induced by the inclusion $\partial\bX\into\bX$ of
the boundary and the space $H_1(\partial_K \bX;\C(\Go)){}\simeq{}\C^2$ is
generated by the meridian~$m$ and Seifert longitude~$\ell$.

\definition[see \cite{DFL2}]\label{def.slope}
If $\Ker i_*$ in~\eqref{eq.boundary}
has dimension one,
it is generated by a single vector $am+b\ell$ for some
$[a:b] \in \mathbb{P}^1(\C)$,
and the \emph{slope} of $K\cup L$ at $\Go\in\CCA(K/L)$ is defined as the quotient
\[*
(K/L)(\Go):=-\frac{a}{b}\in\C\cup\infty.
\]
This notion is extended to all characters $\Go\in\CA(K/L)$ by ``patching'' the
components~$L_i$ on which $\Go_i=1$.
(This operation results in patching with solid tori the corresponding
boundary components of the manifold $X:=S^3\sminus T_L$.)
\enddefinition

\proposition[see \cite{DFL2}]\label{prop.finite.slope}
The slope at a character $\Go\in\CCA(K/L)$ is well defined
if and only if the two inclusion homomorphisms
$H_1(K;\C(\zeta))\to H_1(S^3\sminus L;\C(\zeta))$,
$\zeta=\Go$ or~$\Go\dm$,
are either both trivial or both nontrivial.
The slope is finite, $(K/L)(\Go)\in\C$,
if and only if both homomorphisms are  trivial.
\done
\endproposition

Note also (see \cite[Section~2.4]{DFL2} for details)
that the slope is always defined on a
\emph{unitary} character $\Go\in(S^1)^\mu$: in this case, by twisted
Poincar\'{e} duality, $\Ker \, i_*$ is a Lagrangian subspace of
\[*
H_1(\partial_K \bX;\C(\Go))=H_1(\partial \bX;\C(\Go)),
\]
\cf.~\eqref{eq.boundary},
with respect to
the twisted intersection form and,
hence, $\dim\Ker i_*=1$.

Recall (see, \eg, \cite{Libgober:char.var}) that the
\emph{characteristic varieties} associated with a $\mu$-colored link~$L$ are the
jump loci
\[*
\CV_r(L):=\bigl\{\Go\in(\Cs)^\mu\bigm|\dim H_1(X;\C(\Go))\ge r\bigr\},
 \quad r\ge0.
\]
They are indeed nested algebraic subvarieties:
\[
(\C\units)^\mu=\CV_0\supset\CV_1\supset\CV_2\supset\ldots,
\qquad \text{ with }\CV_1(L)=\bigl\{\Go\bigm|\Delta_L(\Go)=0\bigr\}.
\label{eq.char.var}
\]
The
first \emph{proper} characteristic variety, \ie, the first member $\CV_r$
of the sequence~\eqref{eq.char.var} such that $\CV_r\ne(\C\units)^\mu$,
is denoted by
$\CVmax:=\CVmax(L)$.
This variety depends on~$L$ only and,
if $\Gl:=\vlk(K,L)=0$, it is a proper algebraic subvariety of the torus
$\CA(K/L)=(\C\units)^\mu$ of admissible characters.

\remark \label{remark.linking1}
If
$\lambda:=\vlk(K,L)\ne0$,
the situation is slightly more involved. Let $\lambda=n\lambda'$,
where $\lambda'\in\Z^\mu$ is a primitive vector.
In view
of~\eqref{eq.CCA}, the
variety $\CA(K/L)$
of admissible characters (depending on~$\Gl$ only)
splits over~$\Q$ into
irreducible components
\[*
\CA_d:=\bigl\{\Phi_d(\omega^{\lambda'})=0\bigr\},\quad d\mathrel|n,
\]
where $\Phi_d$ stands for the cyclotomic polynomial, and we should speak
about a separate \emph{first proper characteristic variety}
$\CVmax^{\lambda,d}(L)\subsetneq\CA_d$ for each
component~$\CA_d$. In general,
$\CVmax^{\lambda,d}(L)\ne\CVmax(L)\cap\CA_d$ as $\CVmax(L)$
may contain~$\CA_d$. To keep the notation uniform, we
occasionally extend it to the
case $\lambda=0$ \via\ $\CA_0:=\CA(K/L)$ and $\CVmax^{0,0}(L):=\CVmax(L)$.
\endremark

\theorem[see {\cite[Theorems~3.19 and 3.21]{DFL2}}]\label{th.rational}
Let $\Gl:=\vlk(K,L)$.
For each rational component $\CA_d\subset\CA(K/L)$, the slope restricts to a
rational function, possibly identical~$\infty$, on
the complement
$\CCA_d\sminus\CVmax^{\Gl,d}(L)$.
In other words, the slope gives
rise to an element of the extended function field $\Q(\CA_d)\cup\infty$.

If
$\CVmax^{\Gl,d}(L)=\CV_1(L)\cap\CA_d$,\ie, $\Delta_L$ does not vanish
identically on~$\CA_d$, one has
\[*
(K/L)(\Go)=-\frac{\nabla'(1,\sqrt\Go)}{2\nabla_L(\sqrt\Go)}\in\C\cup\infty,
\]
where $\nabla'$ is the derivative of $\nabla_{K\cup L}(t,\cdot)$ with respect to~$t$.
\done
\endtheorem

\section{Concordance of links} \label{section.conc}

Two oriented $\mu$-colored links
$L^0$ and $L^1$
are \emph{concordant} if there exists
a collection of properly embedded disjoint locally flat cylinders
$A:= A_1\sqcup \ldots \sqcup A_\mu$
in $S^{3}\times [0,1]$ such that
\[*
\partial A_i\cap(S^3\times 0)=-L^{0}_{i} \text{ and } \partial A_i\cap(S^3\times 1)=L^{1}_{i}
\]
for all $i=1,\ldots,\mu$.
(In general, each $A_i$ is a union of
cylinders.)

\subsection{The concordance invariance}
In the study of knot and link concordance, there is a subset of the complex numbers of particular relevance, the
so-called
\emph{Knotennullstellen}. This was first introduced in \cite{NP17} for knots
and extended to the
multi-component
link case in \cite{CNT}. For our purposes, we only need the following
definition. Consider the subset of Laurent polynomials
 \[*
 U := \bigl\{ p \in  \Z[t_1^{\pm 1},\dots, t_\mu^{\pm 1}]  \bigm|
 p(1,\dots,1)= \pm 1 \bigr\}.
 \]
An element $\omega \in \CA(K/L)$
is called a \emph{concordance root} if there is
a polynomial $p \in U$
such that
$p(\omega)=0$. We denote
by $\CA_c(K/L)  \subset \CA(K/L)$
the subset of
admissible characters that
are \emph{not} concordance roots
and abbreviate $\CCA_c(K/L):=\CA_c(K/L)\cap\CCA(K/L)$.
Note that these sets are larger than
the set~$\mathbb T_{!}$ used
in~\cite{CNT},
since we allow for non-unitary characters.

\remark\label{rem.dense}
If $\vlk(K,L)=0$, the set $\CA_c(K/L)$ is dense in
$\CA(K/L)=(\C\units)^\mu$, as it is a countable intersection of Zariski open
sets. In general, $\CA_c(K/L)$ is only dense in the components $\CA_d$ (see
\autoref{remark.linking1}) for which $d$ is a prime power (or $d=1$ as a
special case). Indeed, if $d$ is \emph{not} a prime power, then
$\Phi_d(\cdot)\in U$ and, hence, each point of~$\CA_d$ is a concordance
root.
\endremark

\theorem\label{th.main}
Let $K^0\cup L^0$ and $K^1\cup L^1$ be two concordant $(1,\mu)$-colored links.
Then
$\CA_c(K^0\!/L^0)$ and $\CA_c(K^1\!/L^1)$ coincide as subsets of
$(\C\units)^\mu$
and one has
$$ (K^0\!/L^0)(\Go)= (K^1\!/L^1)(\Go)$$
for any character $\Go \in \CA_c(K^0\!/L^0)$.
\endtheorem

The proof of \autoref{th.main} is postponed till
\autoref{proof.main}.
The next few corollaries are direct consequences of \autoref{th.main}
and \autoref{th.C-complex} below.

\corollary\label{cor.linkingnotzero}
Let $K^0\cup L^0$ and $K^1\cup L^1$ be two concordant $(1,\mu)$-colored links
 such that $\vlk(K^s,L^s)=0$, $s=0,1$.
  Then the slopes
$K^0\!/L^0$ and $K^1\!/L^1$ are equal as elements of
the extended function field $\Q\bigl((\C\units)^\mu\bigr)\cup\infty$.
In particular, $(K^0\!/L^0)(\Go)= (K^1\!/L^1)(\Go)$ for each
character~$\Go$ in
the complement of the \rom(common\rom) first proper characteristic variety
$\CVmax(L^0)=\CVmax(L^1)$.
\done
\endcorollary

\proof
If $L^0$ and $L^1$ are concordant, their nullities coincide
(see \cite[Theorem 7.1]{CF}); hence, so do their first proper characteristic
varieties. Therefore, the statement is an immediate consequence of
\autoref{th.main}, the rationality of the slope given by
\autoref{th.rational}, and
the density of $\CA_c(K/L)$ discussed in \autoref{rem.dense}.
\endproof

\corollary[of \autoref{cor.linkingnotzero} and \autoref{th.rational}]\label{cor.nabla}
Let $K^0\cup L^0$ and $K^1\cup L^1$ be two concordant $(1,\mu)$-colored links
 such that $\vlk(K^s,L^s)=0$ and $\Delta_{L^s}\not\equiv0$, $s=0,1$. Then
\def\qedsymbol{\donesymbol}\pushQED{\qed}
\[*
\frac{\nabla'_{K^0\cup L^0}(1,\bar t)}{\nabla_{L^0}(\bar t)}
 =\frac{\nabla'_{K^1\cup L^1}(1,\bar t)}{\nabla_{L^1}(\bar t)},\quad
 {\bar t}:=(t_1,\ldots,t_\mu).
\qedhere
\]
\endcorollary

\remark\label{rem.lk=/=0}
\latin{A priori}, the conclusions of Corollaries~\ref{cor.linkingnotzero}
or~\ref{cor.nabla} do not need to hold if $\Gl:=\vlk(K^s,L^s)\ne0$: it is not even
obvious that the first proper varieties
$\CVmax^{\Gl,d}(L^s)$
or even their indices
in~\eqref{eq.char.var}
should coincide if $d$ is not a prime power.
(Note though
that we do not know any counterexample, as that would require going far
beyond the known link tables.) The precise statements, based on
Remarks~\ref{remark.linking1} and~\ref{rem.dense} and Theorems~\ref{th.main}
and~\ref{th.rational}, are left to the reader.
\endremark

Recall that a link is \emph{slice} if it is concordant to an unlink. It is a \emph{boundary link} if the components bound a collection of
 mutually disjoint Seifert surfaces in $S^3$.

\corollary\label{cor.slice}
If $K \cup L$ is a slice link,
then $(K/L)(\Go)=0$ for all $\Go$ in $\CA_{c}(K/L)$.
\done
\endcorollary

\corollary \label{cor.boundary}
If $K \cup L$ is concordant to a boundary link,
then $(K/L)(\Go)=0$ for all $\Go$ in $\CA_{c}(K/L)$
\rom(and for any coloring used to define the slope\rom).
\done
\endcorollary

\autoref{cor.boundary} is in fact a particular case of the following statement
(see~\cite{CF} or \autoref{s.C-complex} below for the definition of
$C$-complex).

\corollary \label{cor.disjoint}
If $K \cup L$ is concordant to a $(1,\mu)$-colored link $K' \cup L'$ admitting a
$C$-complex~$F$  for $L$ and a Seifert surface~$S$ for~$K$ disjoint from~$F$, then
$(K/L)(\Go)=0$ for all $\Go \in \CA_c(K/L)$.
\done
\endcorollary

The following example
illustrates that the values of the slope at concordance roots,
that is outside the set $\CA_{c}(K/L)$, might not be invariant
under concordance.
We observe a similar pattern with knot signatures: Knotennullstelle unitary characters are precisely where they fail to be concordance invariants \cite{CL04, NP17}. See \cite{CNT} for the case of colored links.

\example \label{specialval} Let $K \cup L$ be the $(1,1)$-colored
two-component slice link
\href{http://katlas.math.toronto.edu/wiki/L10n36}{L10n36}, where $K$ is the
unknotted component. One has $\nabla_{K \cup L}(t,t_1)=0$ and
$\nabla_L(t_1)=(t_1 -1 + t_1\1)^2$;
hence, by Theorem 3.21 in~\cite{DFL2},
$(K/L)(\Go)=0$ unless $\Go$ is one of the two roots $\Ga_\pm$
of $\nabla_L$,
which agrees with \autoref{th.main} and \autoref{cor.nabla}.
(By definition, $\alpha_\pm \notin \CA_c(K/L)$.)
A computation using Fox
calculus (see \S3.2 in~\cite{DFL2}) gives us $(K/L)(\alpha_\pm) = \infty$.
\endexample

In the proof of \autoref{th.main} we will need the following lemma.
We state it in our more general setting of arbitrary, not necessarily
unitary, characters, but the proof found in~\cite{CNT} extends literally as
it relies on simple homological algebra.

\lemma[Lemma~2.16 in \cite{CNT}]\label{vanish}
Let $k\ge0$ be
an integer.
If $(X,Y)$ is a
\CW-pair
over $B\Z^{\mu}$
such that $H_i(X,Y;\Z)=0$
for all $0 \le i \le k$, then also $H_i(X,Y;\C(\Go))=0$ for
all $0 \le i \le k$
and any character $\Go\in(\C\units)^\mu$ that is not a concordance root.
\done
\endlemma

\subsection{Proof of \autoref{th.main}}\label{proof.main}
To save space, we abbreviate
$H_*^\Go({-}):=H_*({-};\C(\Go))$.

Let  $D \cup A\subset S^3 \times [0,1]$  be the concordance, $\partial D=  -K^0 \sqcup{K^1}$, and consider an open tubular
neighborhood $T_{D \cup A}$ of $D \cup A$ with a fixed trivialisation
extending Seifert framings
(in the tubular neighborhoods
$T_{K^s\cup L^s}:=T_{D\cup A}\cap(S^3\times s)$, $s=0,1$) of the links.
Denote
\[*
U:= S^3 \times [0,1] \sminus T_{A},\qquad
\bU:= S^3 \times [0,1] \sminus T_{D \cup A}
\]
and let
\[*
X^s:= U\cap(S^3\times s),\qquad
\bX^s:= \bU\cap(S^3\times s)
\]
for $s=0,1$.
The inclusions $\bX^s\into\bU$ send
the meridians of $K^s \cup L^s$ to
those
of $D \cup A$.
  The relative Mayer--Vietoris exact sequences
  applied to
\[*
(S^3 \times I, S^3 \times s)
 = (\bU, \bX^s) \cup (\bT_{D \cup A} , \bT_{K^s\cup L^s} )
 = (U, X^s) \cup (\bT_{A} , \bT_{L^s} )
\]
(where $\bT_*$ stands for the closure of a tubular neighborhood $T_*$)
give us
\[
H_*(\bU,\bX^s)=H_*(U,X^s)=0
\label{eq.vanishing}
\]
for $s=0,1$.
In particular, the inclusions $\bX^s \into \bU$ induce isomorphisms
\[
H_1(\bX^0) \stackrel{\simeq}\longrightarrow H_1(\bU)
\stackrel{\simeq} \longleftarrow H_1(\bX^1)
\label{eq.iso}
\]
preserving the meridians and, thus, identify the three character tori.
Since the trivialization of $T_D$ homotopes $\ell^0$ to $\ell^1$,
we have
$\CA_c(K^0\!/L^0)=\CA_c(K^1\!/L^1)$
(\cf. \autoref{rem.charac}).

From now on, patching, if necessary, a few components of both links (and
the concordance), we can assume the character~$\Go$ non-vanishing,
$\Go \in \CCA_c(K^0\!/L^0)$. Referring to \autoref{rem.charac} and using the
above identification of the character tori, we can regard~$\Go$ as a
homomorphism $\pi_1(\bU)\to\C\units$.
The twisted Mayer--Vietoris sequence
applied to the pairs
\[*
(U,X^s) = (\bU , \bX^s) \cup (\bT_D,\bT_{K^s})
\]
gives us, for all $i$,
\[*
\rightarrow H_{i}^\Go(D \times S^1,K^s \times S^1)
\rightarrow H_{i}^\Go(\bU,\bX^s) \oplus H_{i}^\Go(\bT_D,\bT_{K^s})
\rightarrow H_i^\Go(U,X^s) \rightarrow ,
\]
where $\{ \cdot \} \times S^1$ are the meridians of $K^s$ and $D$,
on which $\Go$ is trivial.
Since
\[*
H_{*}^\Go(D \times S^1,K^s \times S^1)=0\quad\text{and}\quad
H_*^\Go(\bU,\bX^s)=H_*^\Go(U,X^s)=0,
\]
the latter by \autoref{vanish} and~\eqref{eq.vanishing},
we obtain
$H_*^\Go(\bU, \bX^s)=0$
and the inclusions $\bX^s \into \bU$ induce isomorphisms
\[*
  H_1^\Go(\bX^0) \stackrel{\simeq} \longrightarrow H_1^\Go(\bU)
  \stackrel{\simeq}  \longleftarrow H_1^\Go(\bX^1)
\]
preserving the meridians
and, similar to~\eqref{eq.iso}, taking the class of~$\ell^0$ to
that of~$\ell^1$.
It follows that $am^0+b\ell^0=0 \in H_1^\Go(\bX^0)$
if and only if $am^1+b\ell^1=0 \in H_1^\Go(\bX^1)$.
\qed

\section{Computation with Seifert forms} \label{section.comp}

In this section, unless specified otherwise, we abbreviate
\[*
H_*({-}):=H_*({-};\C),\quad
H^*({-}):=H^*({-};\C),\quad
H_*^\Go({-})=H_*({-};\C(\Go)).
\]
For a character $\Go\in(\Cs\sminus1)^\mu$, we also abbreviate
$\tGo{i}:=(1-\Go_i\1)$, $1\le i\le\mu$.

\subsection{Seifert forms}\label{s.C-complex}
Let $L=L_1\cup\ldots\cup L_\mu\subset$ be an oriented $\mu$-colored link in $S^3$.
A \emph{$C$-complex} $F$ for $L$
is a collection of Seifert
surfaces~$F_1,\dots,F_\mu$ for the sublinks $L_1,\dots,L_\mu$
that intersect only along (a finite number of)
\emph{clasps}.
Each class in $H_1(F;\Z)$ can be
represented by a collection of \emph{proper loops},
\ie,
loops $\Ga\:S^1\to F$ such that the pull-back of each clasp is a single
segment (possibly empty). We routinely identify classes, loops, and their
images.

Given a vector $\Ge\in\{\pm1\}^\mu$, the \emph{push-off} $\Ga^\Ge$ of a
proper loop~$\Ga$ is the loop in $S^3\sminus F$ obtained by a slight shift
of~$\Ga$ off each surface~$F_i$ in the direction of~$\Ge_i$. (If $\Ga$ runs
along a clasp $\clasp\subset F_i\cap F_j$, the shift respects both
directions~$\Ge_i$ and~$\Ge_j$.) Due to~\cite{CF}, this
operation gives rise to a well-defined homomorphism
\[*
\Theta^\Ge\:H_1(F;\Z)\to H_1(S^3\sminus F;\Z)=H^1(F;\Z)
\]
(we use Alexander duality),
which can be computed by means of
the
\emph{Seifert forms}
\[*
\theta^\Ge\:H_1(F;\Z)\otimes H_1(F;\Z)\to\Z,\qquad
\Ga\otimes\Gb\mapsto\lk(\Ga,\Gb^\Ge).
\]
Now, given a character $\Go\in(\C\units\sminus1)^\mu$, we define
\[*
\Pi(\Go):=\prod_{i=1}^\mu(1-\Go_i)\in\Cs,\quad
A(\Go):=\sum_{\Ge\in\{\pm1\}^\mu}
 \prod_{i=1}^\mu\Ge_i\Go_i^{(1-\Ge_i)/2}\Theta^\Ge
 \:H_1(F)\to H^1(F)
\]
and let
\begin{equation} \label{operatorE}
\E(\Go):=\Pi(\Go\1)\1A(\Go\1)\:H_1(F)\to H^1(F).
\end{equation}
Throughout the text we will use the shortcut notation $\Ker \E(\Go)^\perp$ to denote the subset of $H^{1}(F)$ defined as $\operatorname{Ann}\Ker\E(\Go)$.
It is straightforward that
\[*
\E^*(\Go)=\E(\Go\1),\quad \bar\E(\Go)=\E(\bar\Go),
\]
where: $E^{*}$ is the adjoint in the sense of linear algebra over an arbitrary field, and for a linear map $L\:U\otimes\C\to V\otimes\C$ between two
complexified real vector spaces, we let
$\bar L\:u\mapsto\overline{L(\bar u)}$.
In particular, if $\Go\in(S^1\sminus1)^\mu$ is unitary, the operator
$\E(\Go)$ is Hermitian, \ie, $\bar\E^*(\Go)=\E(\Go)$; thus, it has a
well-defined signature.
Furthermore, if $\Go$ is unitary,
the operator $\E(\Go\1)$ differs from $H(\Go)$ considered
in~\cite{CF} by the positive real constant
$\Pi(\Go)\1\Pi(\bar\Go)\1$; hence, the two have the same signature and
nullity
and $\E$ can be used instead of $H$ in the following theorem.

\theorem[see~\cite{CF}]\label{th.Cimasoni.Florens}
If
$\Go\in(S^1\sminus1)^\mu$ is a non-vanishing unitary character, then one has
$\Gs_L(\Go)=\sign\E(\Go)$ and
$\eta_L(\Go)=\dim\Ker\E(\Go)+b_0(F)-1$.
\endtheorem

In the case of a $1$-colored link~$L$, the $C$-complex reduces to a single
Seifert surface~$F$, so that $\theta:=\theta^+$ and $\Theta:=\Theta^+$ are
the classical Seifert form and operator, respectively. Since, in this case,
we obviously have $\theta^-=\theta^*$ and, hence, $\Theta^-=\Theta^*$, the
operator~$\E$ takes the classical form\[*
\E(\Go\1)=(1-\Go)\1(\Theta-\Go\Theta^*).
\]

\subsection{The statement}

Let $K\cup L$ be a $(1,\mu)$-colored link. Assume that $\lambda$, the linking vector between $K$ and $L$, vanishes and
fix a $C$-complex~$F$ for~$L$ disjoint from~$K$.
By Alexander duality $H_1(S^3\sminus F;\Z)=H^1(F;\Z)$,
there is a well-defined
cohomology class
\[*
\kappa:=[K]\in H^1(F;\Z)\subset H^1(F),\qquad \kappa\:\Ga\mapsto\lk(\Ga, K).
\]

\theorem \label{th.C-complex}
Under the above assumptions, for any character $\Go\in \CCA(K/L)$,
consider the operator
$\E(\Go)\:H_1(F)\to H^1(F)$, see \eqref{operatorE}.
Then
\[*
(K/L)(\Go)=\begin{cases}
 -\<\Ga,\kappa\>, & \mbox{if } \kappa\in\Im\E(\Go)\cap\Ker\E(\Go)^\perp,\\
 \infty, & \mbox{if } \kappa\notin\Im\E(\Go)\cup\Ker\E(\Go)^\perp,\\
 \mbox{\rm undefined}, & \mbox{otherwise},
\end{cases}
\]
where, in the first case, $\Ga\in H_1(F)$ is any class such that
$\E(\Go)(\Ga)=\kappa$.
\endtheorem

\begin{example}
Consider the Whitehead link $K\cup L$ with the C-complex $F$ depicted in
\autoref{f.whitehead}, which is simply a genus one Seifert surface for the
knot $L$. We want to compute the slope $(K/L)(\Go)$ using
\autoref{th.C-complex} and to this end we fix the basis $\{a,b\}$  of
$H_{1}(F)$ and $\{a',b'\}$ of $H_{1}(S^{3}\sminus F)=H^{1}(F)$ which are
illustrated in Figure~\ref{f.whitehead}. With respect to these bases we
have:
$$
\theta^{+}=
\begin{bmatrix}
0&0\\
1&1
\end{bmatrix},\quad
A(\Go)=
\begin{bmatrix}
0&-\Go\\
1&1-\Go
\end{bmatrix},\quad
E(\Go)=
\begin{bmatrix}
0&(1-\Go)\1\\
(1-\Go\1)\1&1
\end{bmatrix}.
$$
It is evident from the figure that $\kappa$ is the same class as $a'$. One can easily compute a class $\Ga\in H_1(F)$  such that $\E(\Go)(\Ga)=\kappa$:
$$
E(\Go)
\begin{bmatrix}
(1-\Go^{-1})(\Go-1)\\
1-\Go
\end{bmatrix}=
\begin{bmatrix}
1\\
0
\end{bmatrix}=
\kappa
$$
Finally, we calculate the slope as $-\langle\alpha,\kappa\rangle$, that is,
$$
(K/L)(\Go)=(1-\Go)(1-\Go^{-1}),
$$
which coincides with previous computations using Fox calculus
(see~\cite{DFL2}).
\figure
\begingroup%
  \makeatletter%
  \providecommand\color[2][]{%
    \errmessage{(Inkscape) Color is used for the text in Inkscape, but the package 'color.sty' is not loaded}%
    \renewcommand\color[2][]{}%
  }%
  \providecommand\transparent[1]{%
    \errmessage{(Inkscape) Transparency is used (non-zero) for the text in Inkscape, but the package 'transparent.sty' is not loaded}%
    \renewcommand\transparent[1]{}%
  }%
  \providecommand\rotatebox[2]{#2}%
  \newcommand*\fsize{\dimexpr\f@size pt\relax}%
  \newcommand*\lineheight[1]{\fontsize{\fsize}{#1\fsize}\selectfont}%
  \ifx\svgwidth\undefined%
    \setlength{\unitlength}{188.59501968bp}%
    \ifx\svgscale\undefined%
      \relax%
    \else%
      \setlength{\unitlength}{\unitlength * \real{\svgscale}}%
    \fi%
  \else%
    \setlength{\unitlength}{\svgwidth}%
  \fi%
  \global\let\svgwidth\undefined%
  \global\let\svgscale\undefined%
  \makeatother%
  \begin{picture}(1,0.54636159)%
    \lineheight{1}%
    \setlength\tabcolsep{0pt}%
    \put(0,0){\includegraphics[width=\unitlength,page=1]{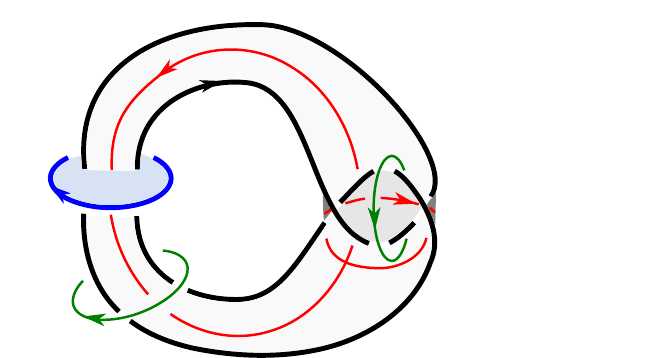}}%
    \put(0.18499041,0.41631695){\color[rgb]{0,0,0}\makebox(0,0)[lt]{\lineheight{1.25}\smash{\begin{tabular}[t]{l}${\scriptstyle\color{red} a}$\end{tabular}}}}%
    \put(-0.00160521,0.26476047){\color[rgb]{0,0,0}\makebox(0,0)[lt]{\lineheight{1.25}\smash{\begin{tabular}[t]{l}$K$\end{tabular}}}}%
    \put(0.05411006,0.06844499){\color[rgb]{0,0.5,0}\makebox(0,0)[lt]{\lineheight{1.25}\smash{\begin{tabular}[t]{l}${\scriptstyle a'}$\end{tabular}}}}%
    \put(0.68492364,0.21587841){\color[rgb]{0,0,0}\makebox(0,0)[lt]{\lineheight{1.25}\smash{\begin{tabular}[t]{l}${\scriptstyle\color{red} b}$\end{tabular}}}}%
    \put(0.56891079,0.33221939){\color[rgb]{0,0.5,0}\makebox(0,0)[lt]{\lineheight{1.25}\smash{\begin{tabular}[t]{l}${\scriptstyle b'}$\end{tabular}}}}%
    \put(0.43263419,0.52978135){\color[rgb]{0,0,0}\makebox(0,0)[lt]{\lineheight{1.25}\smash{\begin{tabular}[t]{l}$L$\end{tabular}}}}%
    \put(0.53524764,0.06044087){\color[rgb]{0,0,0}\makebox(0,0)[lt]{\lineheight{1.25}\smash{\begin{tabular}[t]{l}$F$\end{tabular}}}}%
  \end{picture}%
\endgroup%

\caption{The Whitehead link $K\cup L$ with a C-complex $F$ for $L$ (a Seifert surface in this case) and chosen bases $\{a,b\}$ and $\{a',b'\}$ of $H_{1}(F)$ and $H_{1}(S^{3}\sminus F)=H^{1}(F)$ respectively.}
\label{f.whitehead}
\endfigure
\end{example}

\section{Proof of \autoref{th.C-complex}} \label{section.proof}

\subsection{Geometry of $C$-complexes}\label{s.C-complexes}
The notation and maps introduced in this section are illustrated in \autoref{f.example}.
Let $L$ be a $\mu$-colored link and $F$, a $C$-complex for $L$.
Given a pair $i\ne j$ of indices, let
$C_{ij}:=F_i\cap F_j$ and $\Clasp_{ij}:=\pi_0(C_{ij})$ be the set of clasps in the
intersection of the surfaces $F_i$ and~$F_j$. Denote also $C:=\bigcup C_{ij}$
and
$\Clasp:=\bigcup\Clasp_{ij}$.

 By convention, each clasp $\clasp\in\Clasp_{ij}$ is oriented from
$\clasp\cap L_i$ to $\clasp\cap L_j$, if $i<j$.
The \emph{sign} of $\clasp$, denoted by $\sg\clasp\in\{\pm1\}$, is the local intersection
index $L_i\circ F_j=L_j\circ F_i$ at the corresponding endpoint of $\clasp$.

Fix a regular open neighborhood $V\subset F$ of the union of all clasps,
denote by~$\bV$ its closure, and let
$\Fo_i:=F_i\sminus V$ for all $i$. 
Then, we have
$\partial\Fo_i=\dL\Fo_i\cup\dC\Fo_i$,
where
\[*
\dL\Fo_i:=\partial\Fo_i\cap L,\qquad
\dC\Fo_i:=\partial\Fo_i\cap\bV.
\]
Given a clasp $\clasp\in\Clasp_{ij}$, let $\bV_\clasp$ be the connected
component of~$\bV$ containing~$\clasp$, and let
$\clasp_i\in H_1(\Fo_i,\dL\Fo_i)$ be the arc
$\Fo_i\cap\bV_\clasp$, with its boundary orientation induced
from~$V$, as well as the class realized by this arc.

The following statement is a formalization of the intuitive fact that any
class in $H^1(F)$ can be represented as the intersection index with a certain
surface $S\subset S^3$ such that $\partial S\cap F=\varnothing$;
on the other hand,
any such
surface can be made disjoint from~$C$
and, when doing
so, each clasp can be ``circumvented'' in two ways.

\lemma\label{lem.H1}
The intersection 
pairing establishes an isomorphism
\[*
H^1(F)=\bigoplus_{i=1}^\mu H_1(\Fo_i,\dL\Fo_i)\Big/
 \bigl\{\clasp_i+\clasp_j=0\bigm|\clasp\in\Clasp_{ij},\ 1\le i<j\le\mu\bigr\}.
\]
\endlemma

\proof
Since all groups involved are torsion free,
the statement follows from the exact sequence of the pair $(F,\bV)$:

\[*
0 \longrightarrow H_1(F) \longrightarrow H_1(F,\bV) \longrightarrow H_0(\bV) \longrightarrow H_0(F),
\]
where $H_1(F,\bV)= \bigoplus_i H_1(\Fo_i,\dC\Fo_i)$.
Then, there remains to apply
Poincar\'{e}--Lefschetz duality
$H^1(\Fo_i,\dC\Fo_i)=H_1(\Fo_i,\dL\Fo_i)$.
\endproof

Let $\Ge\in\{\pm1\}^\mu$. Pick a class $\Ga\in H_1(F)$, represent it by a
proper
loop, and denote by $\Ga_i^\Ge\in H_1(\Fo_i,\dL\Fo_i)$ the class
realized by the arc $\Ga\cap F_i$ pushed off each clasp $\clasp\in\Clasp_{ij}$
in the direction
prescribed by~$\Ge_j$.
Passing further to the image in $H^1(F)$, see \autoref{lem.H1},
we obtain a well-defined homomorphism
$\rel_i^\Ge\:H_1(F)\to H^1(F)$. It is easily seen that $\rel_i^\Ge$ is independent of
$\Ge_i$. In fact,
\[*
\rel_i^\Ge\Ga=\Theta^{\Ge[-i]}\Ga-\Theta^{\Ge[+i]}\Ga,
\]
where $\Ge[\pm i]$ is obtained from $\Ge$ by replacing the $i$-th component
by $\pm1$.
\figure
\scalebox{.7}{
\begingroup%
  \makeatletter%
  \providecommand\color[2][]{%
    \errmessage{(Inkscape) Color is used for the text in Inkscape, but the package 'color.sty' is not loaded}%
    \renewcommand\color[2][]{}%
  }%
  \providecommand\transparent[1]{%
    \errmessage{(Inkscape) Transparency is used (non-zero) for the text in Inkscape, but the package 'transparent.sty' is not loaded}%
    \renewcommand\transparent[1]{}%
  }%
  \providecommand\rotatebox[2]{#2}%
  \newcommand*\fsize{\dimexpr\f@size pt\relax}%
  \newcommand*\lineheight[1]{\fontsize{\fsize}{#1\fsize}\selectfont}%
  \ifx\svgwidth\undefined%
    \setlength{\unitlength}{303.67648489bp}%
    \ifx\svgscale\undefined%
      \relax%
    \else%
      \setlength{\unitlength}{\unitlength * \real{\svgscale}}%
    \fi%
  \else%
    \setlength{\unitlength}{\svgwidth}%
  \fi%
  \global\let\svgwidth\undefined%
  \global\let\svgscale\undefined%
  \makeatother%
  \begin{picture}(1,0.68320572)%
    \lineheight{1}%
    \setlength\tabcolsep{0pt}%
    \put(0,0){\includegraphics[width=\unitlength,page=1]{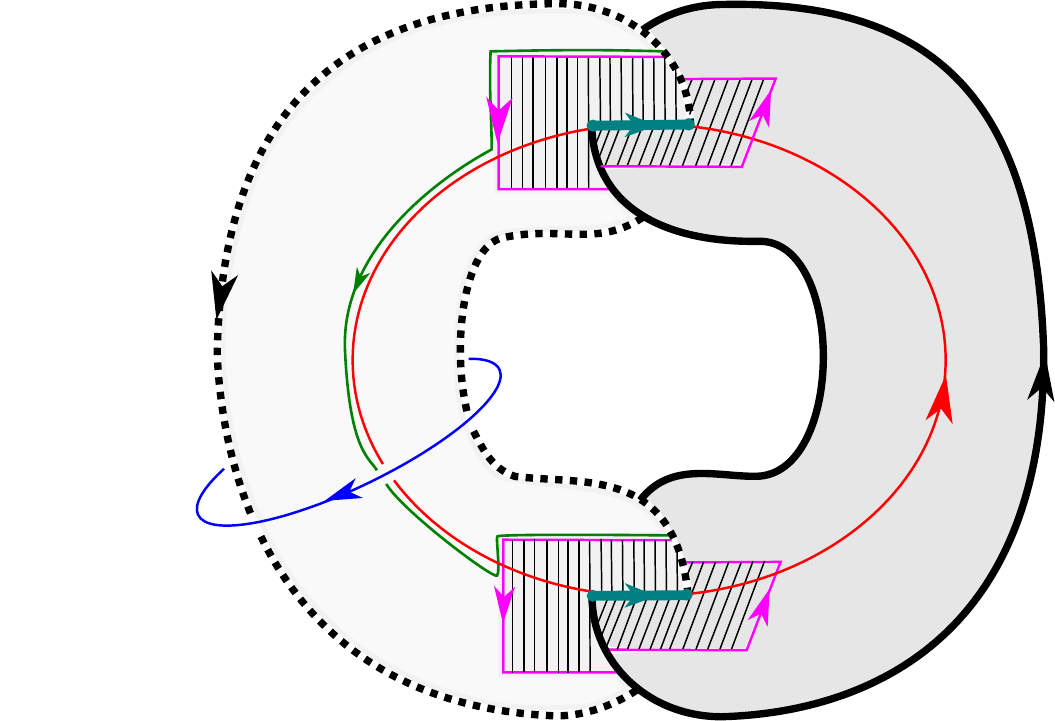}}%
    \put(0.33134341,0.48551254){\color[rgb]{0,0.5,0}\makebox(0,0)[lt]{\lineheight{1.25}\smash{\begin{tabular}[t]{l}$\alpha_1^{++}$\end{tabular}}}}%
    \put(-0.00079276,0.14159095){\color[rgb]{0,0,1}\makebox(0,0)[lt]{\lineheight{1.25}\smash{\begin{tabular}[t]{l}$\Theta^{-+}\alpha-\Theta^{++}\alpha$\end{tabular}}}}%
    \put(0.25720732,0.63410092){\color[rgb]{0,0,0}\makebox(0,0)[lt]{\lineheight{1.25}\smash{\begin{tabular}[t]{l}$L_1$\end{tabular}}}}%
    \put(0.90302319,0.64082453){\color[rgb]{0,0,0}\makebox(0,0)[lt]{\lineheight{1.25}\smash{\begin{tabular}[t]{l}$L_2$\end{tabular}}}}%
    \put(0.35355992,0.08276531){\color[rgb]{0,0,0}\makebox(0,0)[lt]{\lineheight{1.25}\smash{\begin{tabular}[t]{l}$F_1$\end{tabular}}}}%
    \put(0.82761326,0.08456686){\color[rgb]{0,0,0}\makebox(0,0)[lt]{\lineheight{1.25}\smash{\begin{tabular}[t]{l}$F_2$\end{tabular}}}}%
    \put(0.85239572,0.48313781){\color[rgb]{1,0,0}\makebox(0,0)[lt]{\lineheight{1.25}\smash{\begin{tabular}[t]{l}$\alpha$\end{tabular}}}}%
    \put(0.59410349,0.13583064){\color[rgb]{0,0.5,0.5}\makebox(0,0)[lt]{\lineheight{1.25}\smash{\begin{tabular}[t]{l}$\frak b$\end{tabular}}}}%
    \put(0.6028694,0.58067108){\color[rgb]{0,0.5,0.5}\makebox(0,0)[lt]{\lineheight{1.25}\smash{\begin{tabular}[t]{l}$\clasp$\end{tabular}}}}%
    \put(0.48924217,0.64610781){\color[rgb]{0,0,0}\makebox(0,0)[lt]{\lineheight{1.25}\smash{\begin{tabular}[t]{l}$+$\end{tabular}}}}%
    \put(0.68225987,0.64528268){\color[rgb]{0,0,0}\makebox(0,0)[lt]{\lineheight{1.25}\smash{\begin{tabular}[t]{l}$+$\end{tabular}}}}%
    \put(0.43279109,0.56601789){\color[rgb]{1,0,1}\makebox(0,0)[lt]{\lineheight{1.25}\smash{\begin{tabular}[t]{l}$\clasp_1$\end{tabular}}}}%
    \put(0.73573392,0.56398364){\color[rgb]{1,0,1}\makebox(0,0)[lt]{\lineheight{1.25}\smash{\begin{tabular}[t]{l}$\clasp_2$\end{tabular}}}}%
    \put(0.43314396,0.10453056){\color[rgb]{1,0,1}\makebox(0,0)[lt]{\lineheight{1.25}\smash{\begin{tabular}[t]{l}$\frak b_1$\end{tabular}}}}%
    \put(0.7383323,0.09429874){\color[rgb]{1,0,1}\makebox(0,0)[lt]{\lineheight{1.25}\smash{\begin{tabular}[t]{l}$\frak b_2$\end{tabular}}}}%
  \end{picture}%
\endgroup%
}
\caption{This minimal example shows a two colored link $L=L_{1}\cup L_{2}$ bounding a C-complex with two positive claps. In this example $\Clasp=\Clasp_{12}=\{\clasp, \frak b\}$. The lined subset is the open set $V$ with two connected components $V_{\clasp}$ and $V_{\frak b}$. The relative class $\alpha_{1}^{++}\in H_1(\Fo_1,\dL\Fo_1)$ and the element $\Theta^{-+}\alpha-\Theta^{++}\alpha=\rel_{1}^{++}\alpha\in H^{1}(F)$ are identified through the isomorphism in \autoref{lem.H1}.}
\label{f.example}
\endfigure
Furthermore, for an index $j\ne i$, we have
\[
\rel_i^{\Ge[+j]}\Ga-\rel_i^{\Ge[-j]}\Ga
 =\rel_{ij}\Ga
 :=\sum_{\clasp\in\Clasp_{ij}}\sg\clasp\cdot\<\Ga,\clasp_i\>\clasp_i.
\label{eq.rel}
\]
For the reader's convenience a local illustration is presented in \autoref{f.relij}.
(Note that $\<\Ga,\clasp_i\>\clasp_i=\<\Ga,\clasp_j\>\clasp_j$ for each clasp
$\clasp\in\Clasp_{ij}$ and, hence, $\rel_{ij}\Ga=\rel_{ji}\Ga$ as elements of $H^{1}(F)$.)
\figure
\scalebox{.7}{
\begingroup%
  \makeatletter%
  \providecommand\color[2][]{%
    \errmessage{(Inkscape) Color is used for the text in Inkscape, but the package 'color.sty' is not loaded}%
    \renewcommand\color[2][]{}%
  }%
  \providecommand\transparent[1]{%
    \errmessage{(Inkscape) Transparency is used (non-zero) for the text in Inkscape, but the package 'transparent.sty' is not loaded}%
    \renewcommand\transparent[1]{}%
  }%
  \providecommand\rotatebox[2]{#2}%
  \newcommand*\fsize{\dimexpr\f@size pt\relax}%
  \newcommand*\lineheight[1]{\fontsize{\fsize}{#1\fsize}\selectfont}%
  \ifx\svgwidth\undefined%
    \setlength{\unitlength}{291.76536205bp}%
    \ifx\svgscale\undefined%
      \relax%
    \else%
      \setlength{\unitlength}{\unitlength * \real{\svgscale}}%
    \fi%
  \else%
    \setlength{\unitlength}{\svgwidth}%
  \fi%
  \global\let\svgwidth\undefined%
  \global\let\svgscale\undefined%
  \makeatother%
  \begin{picture}(1,0.34937105)%
    \lineheight{1}%
    \setlength\tabcolsep{0pt}%
    \put(0,0){\includegraphics[width=\unitlength,page=1]{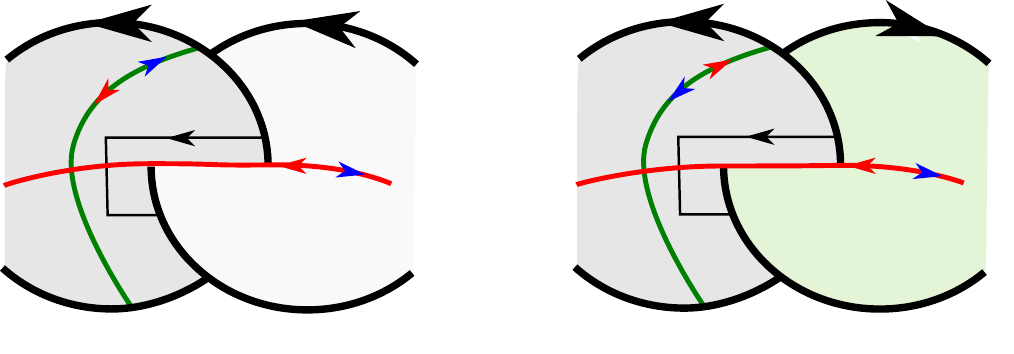}}%
    \put(0.17022956,0.23457377){\color[rgb]{0,0,0}\makebox(0,0)[lt]{\lineheight{1.25}\smash{\begin{tabular}[t]{l}$\clasp_i$\end{tabular}}}}%
    \put(0.00861917,0.10165674){\color[rgb]{0,0,0}\makebox(0,0)[lt]{\lineheight{1.25}\smash{\begin{tabular}[t]{l}$F_i$\end{tabular}}}}%
    \put(0.36782777,0.1025454){\color[rgb]{0,0,0}\makebox(0,0)[lt]{\lineheight{1.25}\smash{\begin{tabular}[t]{l}$F_j$\end{tabular}}}}%
    \put(0.09862991,0.00796603){\color[rgb]{0,0.5,0}\makebox(0,0)[lt]{\lineheight{1.25}\smash{\begin{tabular}[t]{l}$\rel_{ij}\alpha$\end{tabular}}}}%
    \put(0.87226921,0.19636151){\color[rgb]{1,0,0}\makebox(0,0)[lt]{\lineheight{1.25}\smash{\begin{tabular}[t]{l}$\alpha$\end{tabular}}}}%
    \put(0.05736996,0.2870046){\color[rgb]{0,0,0}\makebox(0,0)[lt]{\lineheight{1.25}\smash{\begin{tabular}[t]{l}$+$\end{tabular}}}}%
    \put(0.26887056,0.2896706){\color[rgb]{0,0,0}\makebox(0,0)[lt]{\lineheight{1.25}\smash{\begin{tabular}[t]{l}$+$\end{tabular}}}}%
    \put(0.73936335,0.23593273){\color[rgb]{0,0,0}\makebox(0,0)[lt]{\lineheight{1.25}\smash{\begin{tabular}[t]{l}$\clasp_i$\end{tabular}}}}%
    \put(0.57419834,0.10263409){\color[rgb]{0,0,0}\makebox(0,0)[lt]{\lineheight{1.25}\smash{\begin{tabular}[t]{l}$F_i$\end{tabular}}}}%
    \put(0.93264314,0.09660434){\color[rgb]{0,0,0}\makebox(0,0)[lt]{\lineheight{1.25}\smash{\begin{tabular}[t]{l}$F_j$\end{tabular}}}}%
    \put(0.66065457,0.00577038){\color[rgb]{0,0.5,0}\makebox(0,0)[lt]{\lineheight{1.25}\smash{\begin{tabular}[t]{l}$\rel_{ij}\alpha$\end{tabular}}}}%
    \put(0.30214373,0.20127511){\color[rgb]{1,0,0}\makebox(0,0)[lt]{\lineheight{1.25}\smash{\begin{tabular}[t]{l}$\alpha$\end{tabular}}}}%
    \put(0.60873064,0.28303166){\color[rgb]{0,0,0}\makebox(0,0)[lt]{\lineheight{1.25}\smash{\begin{tabular}[t]{l}$+$\end{tabular}}}}%
    \put(0.90458789,0.28195221){\color[rgb]{0,0,0}\makebox(0,0)[lt]{\lineheight{1.25}\smash{\begin{tabular}[t]{l}$-$\end{tabular}}}}%
  \end{picture}%
\endgroup%
}
\caption{In this figure, the element $\alpha\in H_{1}(F)$ is depicted with both possible orientations. The orientation of the element $\rel_{ij}\alpha$ depends on the sign of the clasp, as illustrated. Remark that the element $\rel_{ij}\alpha$ is by definition in $H^{1}(F)$: the green curve depicted is a representative of that element \via\ \autoref{lem.H1}.}
\label{f.relij}
\endfigure
Let ${-}:=[-1,\ldots,-1]\in\{\pm1\}^\mu$.
Then, applying the last two equations inductively, for each
$\Ge\in\{\pm1\}^\mu$ we get\[
\Theta^\Ge\Ga-\Theta^\mone\Ga
 =-\sum_{i\,:\,\Ge_i>0}\rel_i^\mone\Ga
 -\sum_{i<j\,:\,\Ge_i=\Ge_j>0}\rel_{ij}\Ga.
\label{eq.Theta-1}
\]

\remark\label{rem.rel}
It follows from~\eqref{eq.Theta-1} that, as in the classical case of a single
Seifert surface, all operators~$\Theta^\Ge$ are almost determined by any one
of them, as the relativization homomorphisms
$\rel_i^\Ge$ and $\rel_{ij}$ are intrinsic to the
abstract $C$-complex~$F$ \emph{with prescribed signs $\sg\clasp$} of the
clasps.
In the classical case, \eqref{eq.Theta-1} takes the well-known
form
\[*
\Theta^*-\Theta=\rel\:H_1(F)\to H_1(F,\partial F)=H^1(F),
\]
which explains the notation~$\rel$.
\endremark

Now, given a character $\Go\in(\Cs\sminus1)^\mu$, observe that
\[*
A(\Go)=\Pi(\Go)\Theta^\mone
 +\sum_{\Ge\in\{\pm1\}^\mu}
 \prod_{i=1}^\mu\Ge_i\Go_i^{(1-\Ge_i)/2}(\Theta^\Ge-\Theta^\mone).
\]
Hence, using~\eqref{eq.Theta-1}, rearranging the terms,
and using the definition $\tGo{i}=1-\Go_i\1$,
we arrive at
\[
\E(\Go)=\Theta^\mone-\dif(\Go),\quad
\dif(\Go):=
 \sum_{i=1}^\mu\tGo{i}\1\rel_i^\mone
 +\sum_{1\le i<j\le\mu}\tGo{i}\1\tGo{j}\1\rel_{ij}.
\label{eq.A=rel}
\]

\subsection{Reference sheets}\label{s.ref.sheet}
We
briefly recall how twisted homology can be computed \via\ coverings.
Consider a connected \CW-complex~$X$, an abelian group~$G$,
and an epimorphism $\Gf\:\pi_1(X) \onto H_1(X;\Z)\onto G$.
The kernel of $\Gf$,
which is a normal
  subgroup of $\pi_1(X)$,
gives rise to a Galois
$G$-covering $\tilde X\to X$, where the deck transformation $g{{}\in G}$
sends
a point $\tilde x\in\tilde X$
to the other endpoint of the arc that begins at $\tilde x$ and covers a loop representing
    an element of $\Gf^{-1}(g)$. This
model induces a structure of $\Z[G]$-module
    on $C_*(\tilde X)$ and,
for each multiplicative character $\Go\:G\to\C\units$,
there is
a canonical
chain isomorphism of complexes of $\C(\Go)$-modules
\[*
C_*(X;\C(\Go)) \simeq C_*(\tilde X)\otimes_{\Z G}\C(\Go).
\]
 Occasionally, the
homomorphism $\Gf\:H_1(X;\Z)\to G$ might not necessarily be surjective.
(Typically, this situation occurs when we restrict the construction to a
subcomplex $Y\subset X$.)
Then, letting $G':=\Im\Gf$, the $G$-covering $\tilde X$ consists of $[G:G']$
connected
components, each isomorphic to the $G'$-covering $\tilde X'$, and we have
\[*
C_*(\tilde X)\simeq C_*(\tilde X')\otimes_{\Z G'}\Z G.
\]
However, this isomorphism is
no longer
canonical: to make it such, we need to fix a
\emph{reference component} $\tilde X'\subset\tilde X$. An important special
case is that where the restriction of~$\Go$ to $X$ is trivial. Then we have
an isomorphism
\[*
H_*\bigl(C_*(\tilde X)\otimes_{\Z G}\C(\Go)\bigr)\simeq H^\Go_*(X)=H_*(X),
\]
canonical {\em provided that a \emph{reference sheet} $X$ in
the trivial covering $\tilde X\to X$ is fixed}.

Back to the original set-up,
when dealing with the twisted homology, we need to avoid the ramification
locus~$L$. Hence, we
fix pairwise disjoint tubular neighborhoods
$T_i\supset L_i$ and, denoting by $\bT_i$ the closure of~$T_i$ and
letting $T:=\bigcup_iT_i$, $\bT:=\bigcup_i\bT_i$, introduce
\[*
\gathered
S_L:=S^3\sminus T,\qquad
F_L:=(F\cup\bT)\sminus T\subset S_L,\qquad
C_L:=C\sminus T,\\
\bV_L:=\bV\sminus T,\qquad
\dL\bV_L:=\bV_L\cap\bT,
\endgathered
\]
see \autoref{f.complement}.
\figure
\scalebox{.6}{
\begingroup%
  \makeatletter%
  \providecommand\color[2][]{%
    \errmessage{(Inkscape) Color is used for the text in Inkscape, but the package 'color.sty' is not loaded}%
    \renewcommand\color[2][]{}%
  }%
  \providecommand\transparent[1]{%
    \errmessage{(Inkscape) Transparency is used (non-zero) for the text in Inkscape, but the package 'transparent.sty' is not loaded}%
    \renewcommand\transparent[1]{}%
  }%
  \providecommand\rotatebox[2]{#2}%
  \newcommand*\fsize{\dimexpr\f@size pt\relax}%
  \newcommand*\lineheight[1]{\fontsize{\fsize}{#1\fsize}\selectfont}%
  \ifx\svgwidth\undefined%
    \setlength{\unitlength}{372.75350831bp}%
    \ifx\svgscale\undefined%
      \relax%
    \else%
      \setlength{\unitlength}{\unitlength * \real{\svgscale}}%
    \fi%
  \else%
    \setlength{\unitlength}{\svgwidth}%
  \fi%
  \global\let\svgwidth\undefined%
  \global\let\svgscale\undefined%
  \makeatother%
  \begin{picture}(1,0.5496345)%
    \lineheight{1}%
    \setlength\tabcolsep{0pt}%
    \put(0,0){\includegraphics[width=\unitlength,page=1]{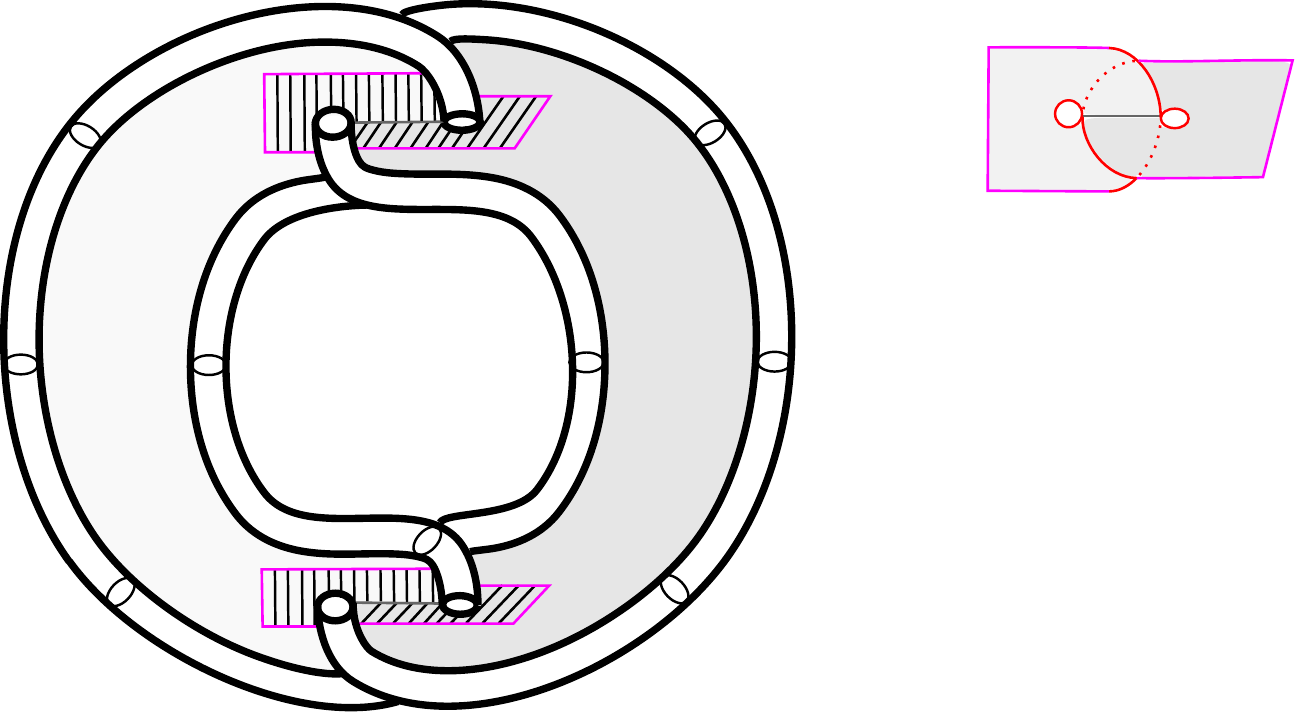}}%
  \end{picture}%
\endgroup%
}
\caption{A minimal example of the set $F_{L}=(F\cup\bT)\sminus T$
consisting of the gray shaded surface together with the two depicted tori. The lined subset is $\bar V_{L}$. To the right we have a copy of a connected component of $\bar V_{L}$ with the subset $\partial_{L}\bar V_{L}$ highlighted in red.}
\label{f.complement}
\endfigure
Here, $V\supset C$ is the neighborhood introduced in \autoref{s.C-complexes},
and we assume the radius of~$T$ so small that $F_{i}\cap\bT_{j}\subset V$ for
each $i\ne j$.

Formally, we also need to shrink the surfaces $\Fo_i$ to
$\Fo_i\sminus T$, changing the boundary $\dL\Fo_i$ to
$(\Fo_i\sminus T)\cap\bT$; however, using the obvious isomorphisms
in (co-)homology, we keep the notation $(\Fo_i,\dL\Fo_i)$ for
these new pairs.

We make use of the isomorphisms
\begin{align}
H^\Go_*(S_L,F_L)&\simeq H_*(S_L,F_L)=H_*(S,F),
\label{eq.(S,F)}\\
H^\Go_*(\Fo_i,\dL\Fo_i)&\simeq H_*(\Fo_i,\dL\Fo_i),
\label{eq.(F,dF)}\\
H^\Go_*(\bV_L,\dL\bV_L)=H^\Go_*(C_L,\partial C_L)&\simeq H_*(C_L,\partial C_L)
 =H_*(C,\partial C),
\label{eq,(V,dV)}
\end{align}
\etc.,
and, in order to fix the (not quite canonical in the context of
a common $G$-covering) isomorphisms denoted by
$\simeq$, we need a coherent choice of reference sheets, upon which we change
the notation to $=$.
(The other isomorphisms are standard combinations of excision and homotopy
equivalences and, thus, are canonical.)
To this end, we consider a ``negative''
collar
(trace of the push-off in the negative direction)
$N:=(-\Gd,0)\times(F\sminus T)$, $\Gd\ll\operatorname{radius}(\bT)$, and,
letting $S'_L:=S_L\sminus N$,
use excision to identify
\[*
H_*(S_L,F_L)=H_*(S'_L,\partial S'_L),\qquad
H^\Go_*(S_L,F_L)=H^\Go_*(S'_L,\partial S'_L).
\]
Since the covering is obviously trivial over
$S'_L$, we can choose and fix a reference sheet
$S'_L\subset\tilde S_L$ and use it for~\eqref{eq.(S,F)}.
There remains to observe that this sheet contains a single copy of each
of~$\Fo_i$ and $C_L$, which are used for
\eqref{eq.(F,dF)} and~\eqref{eq,(V,dV)}, respectively.

\convention\label{conv.partial}
We have then $H_2^\Go(S_L,F_L)=H_2(S_L,F_L)$ and
$H_1(F_L)=H_1^\Go(F_L)$. For the twisted boundary operators like
\[*
H_2(S_L,F_L)\to H_1(F_L),
\]
we assume that $\partial^\Go=\sum_i(\partial^-+\Go_i\1\partial^+)$,
where $\partial^{+}$ is the \emph{lower} boundary
(the $+$ superscript is related to the orientation conventions.)
\endconvention

\convention\label{conv.loops}
The ``reference lift'' of a loop is the loop in the covering whose \emph{end point} is in the reference sheet.
\endconvention

\subsection{The homology of $F$}\label{s.twisted}
Throughout this section, \emph{we assume that $F$ is connected and that $\kappa\ne0$}.
Recall from \autoref{lem.H1} that $H^1(F)$ is a quotient of $\bigoplus H_1(\Fo_i,\dL\Fo_i)$ by relations of the form
$\clasp_i+\clasp_j=0$.
We deduce the following description of the twisted homology of $F$.

\lemma\label{lem.tilde.H1}
The assignment $\tau\:H^1(F)\to H_1^\Go(F_L,\partial\bT)=H_1^\Go(F_L)$
\[*
\sum_{i=1}^\mu\Ga_i\longmapsto
 \operatorname{inclusion}_*\bigoplus_{i=1}^\mu\tGo{i}\Ga_i,\quad
 \Ga_i\in H_1(\Fo_i,\dL\Fo_i),
\]
is a well-defined isomorphism.
\endlemma

\proof
The isomorphisms $H_*^\Go(F_L,\partial\bT)=H_*^\Go(F_L)$ follow from
the assumption
$\Go_i\ne1$ for each~$i$ and, hence, $H_*^\Go(\partial\bT)=0$.
We compute $H_1^\Go(F_L,\partial\bT)$ using the relative
Mayer--Vietoris sequence associated to the decomposition
$F\sminus T=\bV_L\cup\bigl(\bigcup_{i=1}^\mu\Fo_i\bigr)$:
\[
H_1^\Go(\partial \bV_L,\dL \bV_L)\longto
 H_1^\Go(\bV_L,\dL \bV_L)\oplus
  \bigoplus_{i=1}^\mu H_1^\Go(\Fo_i,\dL\Fo_i)\overset{p}\longto
 H_1^\Go(F_L,\partial\bT)\to0,
\label{eq.M-V}
\]
the last term
being $H_0^\Go(\partial\bV_L,\dL\bV_L)=0$, see~\eqref{eq,(V,dV)}
and \autoref{f.complement}.
By~\eqref{eq,(V,dV)}, we also have
$H_1^\Go(\partial\bV_L,\dL\bV_L)=\bigoplus\C\clasp_i$,
the summation running over all
$\clasp\in\Clasp_{ij}$ and all pairs $1\le i\ne j\le\mu$.
The inclusions induce the homomorphisms
\[
\alignedat3
\clasp_i&\mapsto\clasp_i&&\in H_1^\Go(\Fo_i,\dL\Fo_i)
 &&=H_1(\Fo_i,\dL\Fo_i),\ \text see~\eqref{eq.(F,dF)},\\
\clasp_i&\mapsto\sg(j-i)\cdot\sg\clasp\cdot\tGo{j}\clasp
 &&\in H_1^\Go(\bV_L,\dL\bV_L)&&={\textstyle\bigoplus_{\clasp\in\Clasp}\C\clasp}.
\endalignedat
\label{eq.inclusion}
\]
(To follow the above formulas, the reader might find helpful the schematics of the behavior of the twisted homology in Figure~\ref{f.schematics}.)
\figure
\scalebox{.7}{
\begingroup%
  \makeatletter%
  \providecommand\color[2][]{%
    \errmessage{(Inkscape) Color is used for the text in Inkscape, but the package 'color.sty' is not loaded}%
    \renewcommand\color[2][]{}%
  }%
  \providecommand\transparent[1]{%
    \errmessage{(Inkscape) Transparency is used (non-zero) for the text in Inkscape, but the package 'transparent.sty' is not loaded}%
    \renewcommand\transparent[1]{}%
  }%
  \providecommand\rotatebox[2]{#2}%
  \newcommand*\fsize{\dimexpr\f@size pt\relax}%
  \newcommand*\lineheight[1]{\fontsize{\fsize}{#1\fsize}\selectfont}%
  \ifx\svgwidth\undefined%
    \setlength{\unitlength}{426.04834586bp}%
    \ifx\svgscale\undefined%
      \relax%
    \else%
      \setlength{\unitlength}{\unitlength * \real{\svgscale}}%
    \fi%
  \else%
    \setlength{\unitlength}{\svgwidth}%
  \fi%
  \global\let\svgwidth\undefined%
  \global\let\svgscale\undefined%
  \makeatother%
  \begin{picture}(1,0.42423951)%
    \lineheight{1}%
    \setlength\tabcolsep{0pt}%
    \put(0,0){\includegraphics[width=\unitlength,page=1]{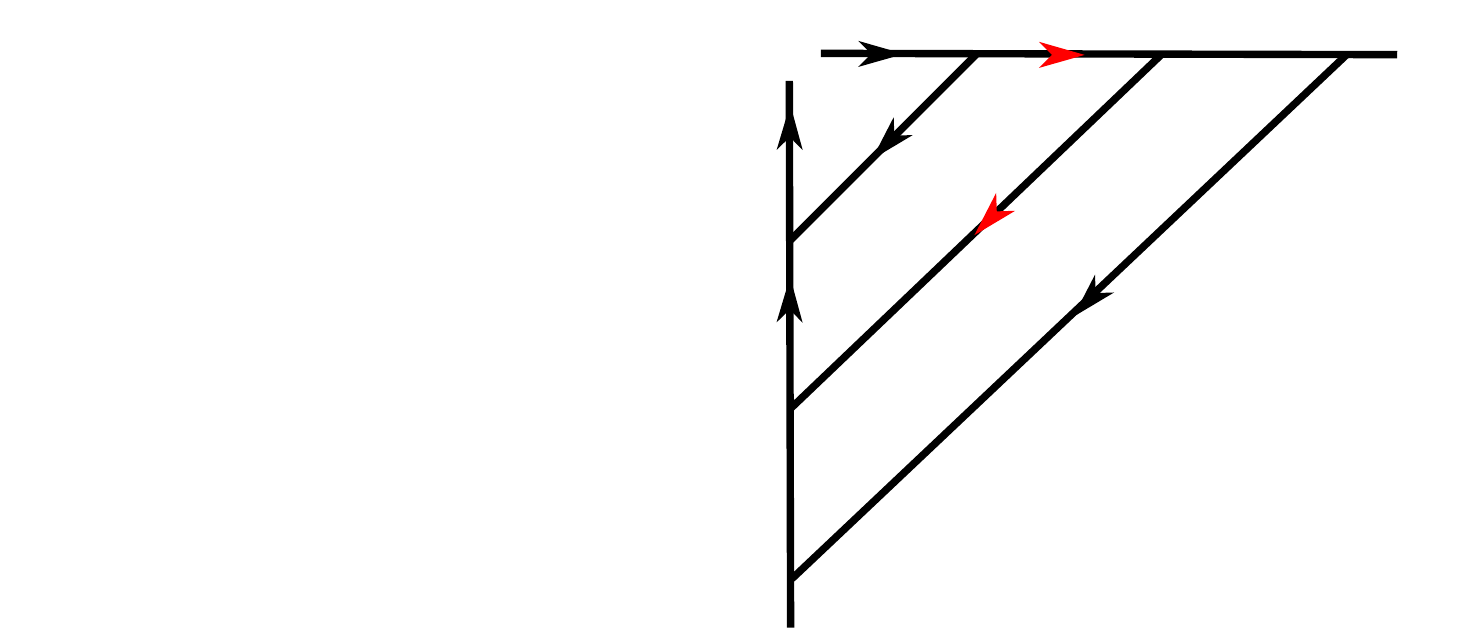}}%
    \put(0.70305869,0.40864589){\color[rgb]{0,0,0}\makebox(0,0)[lt]{\lineheight{1.25}\smash{\begin{tabular}[t]{l}$m_j$\end{tabular}}}}%
    \put(0.81243114,0.40843958){\color[rgb]{0,0,0}\makebox(0,0)[lt]{\lineheight{1.25}\smash{\begin{tabular}[t]{l}$\omega_jm_j$\end{tabular}}}}%
    \put(0.48994797,0.0963286){\color[rgb]{0,0,0}\makebox(0,0)[lt]{\lineheight{1.25}\smash{\begin{tabular}[t]{l}$m_i$\end{tabular}}}}%
    \put(0.46382566,0.21293309){\color[rgb]{0,0,0}\makebox(0,0)[lt]{\lineheight{1.25}\smash{\begin{tabular}[t]{l}$\omega_im_i$\end{tabular}}}}%
    \put(0.46249608,0.32788878){\color[rgb]{0,0,0}\makebox(0,0)[lt]{\lineheight{1.25}\smash{\begin{tabular}[t]{l}$\omega_i^{2}m_i$\end{tabular}}}}%
    \put(0,0){\includegraphics[width=\unitlength,page=2]{schematics.pdf}}%
    \put(0.67910672,0.25952494){\color[rgb]{0,0,0}\makebox(0,0)[lt]{\lineheight{1.25}\smash{\begin{tabular}[t]{l}$\clasp$\end{tabular}}}}%
    \put(0.63506436,0.34138687){\color[rgb]{0,0,0}\makebox(0,0)[lt]{\lineheight{1.25}\smash{\begin{tabular}[t]{l}$\omega_j^{-1}\clasp$\end{tabular}}}}%
    \put(0,0){\includegraphics[width=\unitlength,page=3]{schematics.pdf}}%
    \put(0.56723445,0.27013941){\color[rgb]{0,0,0}\makebox(0,0)[lt]{\lineheight{1.25}\smash{\begin{tabular}[t]{l}$\omega_i\clasp$\end{tabular}}}}%
    \put(0.60001503,0.07480219){\color[rgb]{0,0,0}\makebox(0,0)[lt]{\lineheight{1.25}\smash{\begin{tabular}[t]{l}$\omega_i^{-1}\clasp$\end{tabular}}}}%
    \put(0.85360519,0.30820856){\color[rgb]{0,0,0}\makebox(0,0)[lt]{\lineheight{1.25}\smash{\begin{tabular}[t]{l}$\omega_j\clasp$\end{tabular}}}}%
    \put(0.37187109,0.19657851){\color[rgb]{0,0,0}\makebox(0,0)[lt]{\lineheight{1.25}\smash{\begin{tabular}[t]{l}$\clasp_j$\end{tabular}}}}%
    \put(0,0){\includegraphics[width=\unitlength,page=4]{schematics.pdf}}%
    \put(-0.00066931,0.20458052){\color[rgb]{0,0,0}\makebox(0,0)[lt]{\lineheight{1.25}\smash{\begin{tabular}[t]{l}$\clasp_i$\end{tabular}}}}%
    \put(0.17181954,0.23303233){\color[rgb]{0,0,0}\makebox(0,0)[lt]{\lineheight{1.25}\smash{\begin{tabular}[t]{l}$\clasp$\end{tabular}}}}%
    \put(0.23850337,0.18413082){\color[rgb]{0,0,0}\makebox(0,0)[lt]{\lineheight{1.25}\smash{\begin{tabular}[t]{l}$m_i$\end{tabular}}}}%
    \put(0.10524296,0.18676245){\color[rgb]{0,0,0}\makebox(0,0)[lt]{\lineheight{1.25}\smash{\begin{tabular}[t]{l}$m_j$\end{tabular}}}}%
    \put(0,0){\includegraphics[width=\unitlength,page=5]{schematics.pdf}}%
    \put(0.03578451,0.27393176){\color[rgb]{0,0,0}\makebox(0,0)[lt]{\lineheight{1.25}\smash{\begin{tabular}[t]{l}$F_i$\end{tabular}}}}%
    \put(0.34146237,0.25444265){\color[rgb]{0,0,0}\makebox(0,0)[lt]{\lineheight{1.25}\smash{\begin{tabular}[t]{l}$F_j$\end{tabular}}}}%
    \put(0.55494121,0.4095606){\color[rgb]{0,0,0}\makebox(0,0)[lt]{\lineheight{1.25}\smash{\begin{tabular}[t]{l}$\omega_j^{-1}m_j$\end{tabular}}}}%
    \put(0,0){\includegraphics[width=\unitlength,page=6]{schematics.pdf}}%
  \end{picture}%
\endgroup%
}
\caption{To the left is a local picture of a positive clasp with $i<j$.
To the right, the schematics of the behavior of the lifted curves on a covering space.
Shown in red are the chosen reference lifts.}
\label{f.schematics}
\endfigure
Identifying the two images of each generator~$\clasp_i$, we conclude
that the inclusions $\Fo_i\into F_L$ induce an isomorphism
\[*
\bigoplus_{i=1}^\mu H_1(\Fo_i,\dL\Fo_i)\big/
 \bigl\{\tGo{i}\clasp_i+\tGo{j}\clasp_j=0\bigm|\clasp\in\Clasp_{ij}\bigr\}
 =H_1^\Go(F_L,\partial\bT),
\]
and the isomorphism in the statement follows from
\autoref{lem.H1}.
\endproof

\corollary[of the proof]\label{cor.cylinder}
Given a proper loop $\Ga\subset F$, consider its push-off $\Ga^\mone$ and
its ``trace'' $S^\mone\subset S^3$, \ie,
a
cylinder contained in a regular neighborhood of~$\Ga$ and
such that $S^\mone\cap F=\Ga$ and
$\partial S^\mone=\Ga-\Ga^\mone$. Then, the twisted boundary
$\partial^\Go S^\mone+\Ga^\mone$ is equal to
$\tau\bigl(\dif(\Go)(\Ga)\bigl)\in H_1^\Go(F_L)$,
see~\eqref{eq.A=rel} and \autoref{lem.tilde.H1}.
\endcorollary
\proof
Clearly, $\partial^\Go S^\mone+\Ga^\mone$ is homologous to
the image under~$p$ in~\eqref{eq.M-V} of the cycle
\figure
\scalebox{.7}{
\begingroup%
  \makeatletter%
  \providecommand\color[2][]{%
    \errmessage{(Inkscape) Color is used for the text in Inkscape, but the package 'color.sty' is not loaded}%
    \renewcommand\color[2][]{}%
  }%
  \providecommand\transparent[1]{%
    \errmessage{(Inkscape) Transparency is used (non-zero) for the text in Inkscape, but the package 'transparent.sty' is not loaded}%
    \renewcommand\transparent[1]{}%
  }%
  \providecommand\rotatebox[2]{#2}%
  \newcommand*\fsize{\dimexpr\f@size pt\relax}%
  \newcommand*\lineheight[1]{\fontsize{\fsize}{#1\fsize}\selectfont}%
  \ifx\svgwidth\undefined%
    \setlength{\unitlength}{249.86241569bp}%
    \ifx\svgscale\undefined%
      \relax%
    \else%
      \setlength{\unitlength}{\unitlength * \real{\svgscale}}%
    \fi%
  \else%
    \setlength{\unitlength}{\svgwidth}%
  \fi%
  \global\let\svgwidth\undefined%
  \global\let\svgscale\undefined%
  \makeatother%
  \begin{picture}(1,0.83036142)%
    \lineheight{1}%
    \setlength\tabcolsep{0pt}%
    \put(0,0){\includegraphics[width=\unitlength,page=1]{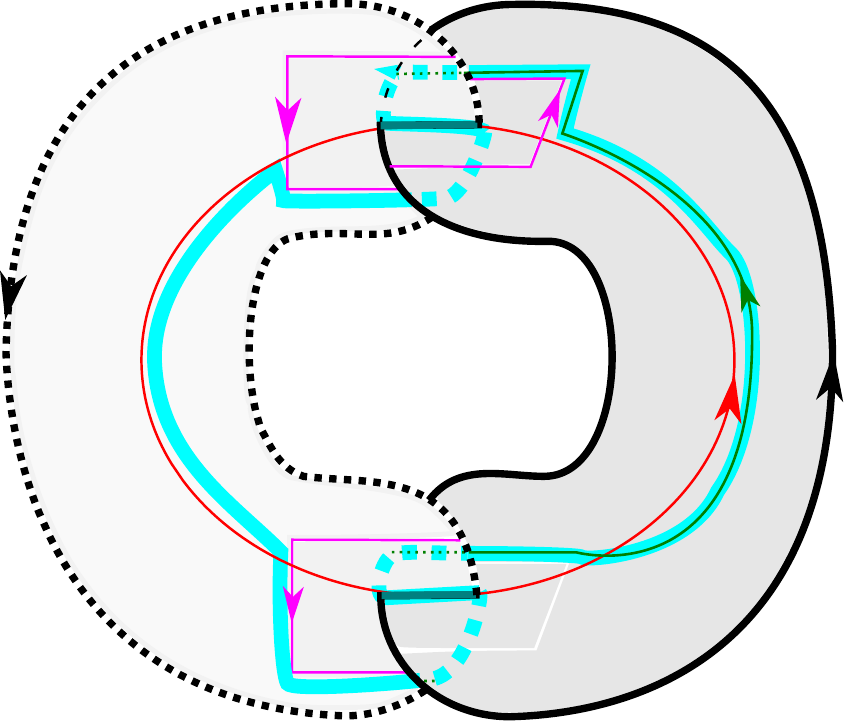}}%
    \put(0.30787094,0.23615344){\color[rgb]{0,0.5,0}\makebox(0,0)[lt]{\lineheight{1.25}\smash{\begin{tabular}[t]{l}$\rel_1^{-}\alpha$\end{tabular}}}}%
    \put(0.06916107,0.77067034){\color[rgb]{0,0,0}\makebox(0,0)[lt]{\lineheight{1.25}\smash{\begin{tabular}[t]{l}$L_1$\end{tabular}}}}%
    \put(0.85406955,0.77884191){\color[rgb]{0,0,0}\makebox(0,0)[lt]{\lineheight{1.25}\smash{\begin{tabular}[t]{l}$L_2$\end{tabular}}}}%
    \put(0.18626564,0.10059077){\color[rgb]{0,0,0}\makebox(0,0)[lt]{\lineheight{1.25}\smash{\begin{tabular}[t]{l}$F_1$\end{tabular}}}}%
    \put(0.76241813,0.10278042){\color[rgb]{0,0,0}\makebox(0,0)[lt]{\lineheight{1.25}\smash{\begin{tabular}[t]{l}$F_2$\end{tabular}}}}%
    \put(0.76986032,0.32752975){\color[rgb]{1,0,0}\makebox(0,0)[lt]{\lineheight{1.25}\smash{\begin{tabular}[t]{l}$\alpha$\end{tabular}}}}%
    \put(0.38914892,0.71325239){\color[rgb]{0,0.5,0.5}\makebox(0,0)[lt]{\lineheight{1.25}\smash{\begin{tabular}[t]{l}$\langle\alpha,\clasp_1\rangle\clasp$\end{tabular}}}}%
    \put(0.35117035,0.78526319){\color[rgb]{0,0,0}\makebox(0,0)[lt]{\lineheight{1.25}\smash{\begin{tabular}[t]{l}$+$\end{tabular}}}}%
    \put(0.58575928,0.7842602){\color[rgb]{0,0,0}\makebox(0,0)[lt]{\lineheight{1.25}\smash{\begin{tabular}[t]{l}$+$\end{tabular}}}}%
    \put(0.2825612,0.68792387){\color[rgb]{1,0,1}\makebox(0,0)[lt]{\lineheight{1.25}\smash{\begin{tabular}[t]{l}$\clasp_1$\end{tabular}}}}%
    \put(0.58214089,0.69745816){\color[rgb]{1,0,1}\makebox(0,0)[lt]{\lineheight{1.25}\smash{\begin{tabular}[t]{l}$\clasp_2$\end{tabular}}}}%
    \put(0.35502962,0.11503718){\color[rgb]{1,0,1}\makebox(0,0)[lt]{\lineheight{1.25}\smash{\begin{tabular}[t]{l}$\frak b_1$\end{tabular}}}}%
    \put(0.64790499,0.1026016){\color[rgb]{1,0,1}\makebox(0,0)[lt]{\lineheight{1.25}\smash{\begin{tabular}[t]{l}$\frak b_2$\end{tabular}}}}%
    \put(0,0){\includegraphics[width=\unitlength,page=2]{zigzag.pdf}}%
    \put(0.79005063,0.62318809){\color[rgb]{0,0.5,0}\makebox(0,0)[lt]{\lineheight{1.25}\smash{\begin{tabular}[t]{l}$\rel_2^{-}\alpha$\end{tabular}}}}%
    \put(0,0){\includegraphics[width=\unitlength,page=3]{zigzag.pdf}}%
    \put(0.37671767,0.47572132){\color[rgb]{1,0,0}\makebox(0,0)[lt]{\lineheight{1.25}\smash{\begin{tabular}[t]{l}$\alpha^-$\end{tabular}}}}%
    \put(0,0){\includegraphics[width=\unitlength,page=4]{zigzag.pdf}}%
    \put(0.45799211,0.10195954){\color[rgb]{0,0.5,0.5}\makebox(0,0)[lt]{\lineheight{1.25}\smash{\begin{tabular}[t]{l}$\langle\alpha,\frak b_1\rangle\frak b$\end{tabular}}}}%
    \put(0,0){\includegraphics[width=\unitlength,page=5]{zigzag.pdf}}%
  \end{picture}%
\endgroup%
}
\caption{The push off $\alpha^{-}$ is to be thought of as located ``behind''
the surface $F_{1}\cup F_{2}$. With the orientations depicted,
together $\alpha$ and $-\alpha^{-}$ are the obvious boundary of the cylinder
$S^{-}$ (not in the picture). The different elements of the cycle described
at the beginning of the proof of Corollary~\ref{cor.cylinder}, $\rel_{i}^{-}\alpha$ and $\langle\alpha,\clasp_{i}\rangle\clasp$, are highlighted.}
\label{f.zigzag}
\endfigure
\[*
\sum_{i=1}^\mu\rel_i^\mone\Ga
  +\sum_{1\le i<j\le\mu}\sum_{\clasp\in\Clasp_{ij}}\<\Ga,\clasp_i\>\clasp,
\]
(see Figure~\ref{f.zigzag} for a simple example.)
Then, by~\eqref{eq.inclusion}, for
all $i<j$ and $\clasp\in\Clasp_{ij}$, we have
$\clasp=\sg\clasp\cdot\tGo{j}\1\clasp_i$ in $H_1^\Go(F_L)$ and,
using \eqref{eq.rel}, we obtain
\begin{multline*}
\sum_{i=1}^\mu\rel_i^\mone\Ga
 +\sum_{1\le i<j\le\mu}\tGo{j}\1
 \sum_{\clasp\in\Clasp_{ij}}\sg\clasp\<\Ga,\clasp_i\>\clasp_i
\overset{\eqref{eq.rel}}{=\joinrel=\joinrel=}
\sum_{i=1}^\mu\rel_i^\mone\Ga+\sum_{1\le i<j\le\mu}\tGo{j}\1\rel_{ij}\Ga\\
=
\sum_{i=1}^\mu\tGo{i}\Biggl(\underbrace{\tGo{i}\1\rel^-_i\Ga
 +\sum_{j=i+1}^\mu\tGo{i}\1\tGo{j}\1\rel_{ij}\Ga}_{R_i}\Biggr).
\end{multline*}
Now, by~\eqref{eq.A=rel}, each $R_i$ is the $i$-th component
of (a representative of) $R(\Go)(\Ga)$, and the statement follows from the
definition of $\tau$ in \autoref{lem.tilde.H1}.
\endproof

We proceed with the computation of the twisted homology of $S_L$ and
$S_L\sminus K$.
We have fixed isomorphisms
\[*
H_*^\Go(S_L,F_L)=H_*(S,F),\qquad
H_*^\Go(S_L\sminus K,F_L)=H_*(S\sminus K,F),
\]
see~\eqref{eq.(S,F)}.
In particular,
\[*
H_1^\Go(S_L,F_L)=H_1^\Go(S_L\sminus K,F_L)=0
\]
(recall that we assume $F$ connected and $\kappa\ne0$) and, by the respective
exact sequences of pairs $(S,F)$ and $(S\sminus K,F)$,
\[*
H_2^\Go(S_L,F_L)=H_1(F),\qquad
H_2^\Go(S_L\sminus K,F_L)=\Ker\kappa\subset H_1(F).
\]
Now, from the
corresponding
twisted exact sequences,
and with the isomorphism~$\tau$ given by \autoref{lem.tilde.H1} taken into
account, we arrive at
\[
H_1^\Go(S_L)=H^1(F)/\Im d,\qquad
H_1^\Go(S_L\sminus K)=H^1(F)/d(\Ker\kappa),
\label{eq.H.twisted}
\]
where $d$ is the composed map
\[
d\:H_1(F)\overset{\partial\1}\longto H_2(S,F)
 =H_2^\Go(S_L,F_L)\overset{\partial^\Go}\longto H_1^\Go(F_L)
 \overset{\tau\1}\longto H^1(F).
\label{eq.d}
\]

\subsection{The twisted homomorphisms}\label{s.homo}
We still assume that $F$ is connected and $\kappa\ne0$.
By~\eqref{eq.H.twisted},
for $X:=S_L$ or $X:=S_L\sminus K$, we have natural
epimorphisms
\[
\pi_X\:H^1(F)\longonto H_1^\Go(X).
\label{eq.pi}
\]
Composing the inclusion with Alexander duality, we obtain a
homomorphism
\[*
\D\:H_1^\Go(X\sminus F_L)=H_1(X\sminus F_L)\to
 H_1(S^3\sminus F)\overset\simeq\longto H^1(F).
\]
Consider also the ``orthogonal projection''
\[*
\alignedat2
\pr_X\:H_1^\Go(X\sminus F_L)&\longto H_1^\Go(X\sminus F_L),\\
\Ga&\longmapsto\Ga&\qquad&\mbox{if $X=S_L$},\\
\Ga&\longmapsto\Ga-\lk(\Ga,K)m&\qquad&\mbox{if $X=S_L\sminus K$}.\\
\endalignedat
\]

\lemma\label{lem.H1.inclusion}
For $X=S_L$ or $S_L\sminus K$ and
any class $\Ga\in H_1^\Go(X\sminus F_L)$, the image of $\pr_X(\Ga)$ under
the inclusion homomorphism
$H_1^\Go(X\sminus F_L)\to H_1^\Go(X)$
is $\pi_X(\D(\Ga))$.
\endlemma
\proof
The statement is a geometric version of \autoref{lem.tilde.H1}. The class
$\Ga':=\pr_X(\Ga)$ is represented by a cycle in $X\sminus F_L$, which bounds a
Seifert surface $G\subset S^{3}\sminus K$.
(This is why we subtract $\lk(\Ga,K)m$ in the case $X=S_L\sminus K$: we want
a Seifert surface disjoint from~$K$.)
Set
$G_{L}:=G\cap S_{L}$.
We can choose the surface $G_{L}$ so that it cuts on $F$ a collection
of arcs $\Ga_i\subset\Fo_i$ with $\partial\Ga_i\subset\dL\Fo_i$.
Then, $\D(\Ga')$ is represented by
\[*
\sum_{i=1}^\mu\Ga_i\in\bigoplus_{i=1}^\mu H_1(\Fo_i,\dL\Fo_i)\longto H^1(F),
\]
see \autoref{lem.H1},
whereas the twisted boundary is
\[
\partial^\Go G_{L}-\Ga'=-\sum_{i=1}^\mu\tGo{i}\Ga_i=-\tau(\D(\Ga)),
\label{eq.dS}
\]
\cf. \autoref{lem.tilde.H1}, implying that $\Ga'=\tau(\D(\Ga))$ in
$H_1^\Go(X)$.
\endproof

\corollary\label{cor.image.K}
For $X=S_L$ or $S_L\sminus K$, let $\Ga\in H_1^\Go(X\sminus F_L)$
be the class of~$[K]$ or $\ell$, respectively. Then, the image of~$\Ga$
in $H_1^\Go(X)$ is $\pi_X(\kappa)$.
\done
\endcorollary

\lemma\label{lem.d}
The homomorphism~$d$ in~\eqref{eq.d} equals $-\E(\Go)$.
\endlemma

\lemma\label{lem.image.m}
For each $\Ga\in H_1(F)$, one has
\[*
\pi_{S_L\sminus K}\bigl(\E(\Go)(\Ga)\bigr)=-\<\Ga,\kappa\>m
\]
in $H_1^\Go(S_L\sminus K)$, see~\eqref{eq.pi}.
\endlemma

\proof[Proof of Lemmas~\ref{lem.d} and~\ref{lem.image.m}]
Let $\Ga\subset F$ be a proper loop and consider its push-off
$\Ga^\mone\subset S^3\sminus(K\cup F)$. Let $S^\mone$ be the trace cylinder as in
\autoref{cor.cylinder}, and let $G$ be a Seifert surface bounded by
$\Ga^\mone$. (For \autoref{lem.image.m}, we replace $\Ga^\mone$ with
its projection
$\pr(\Ga^\mone)=\Ga^\mone-\<\Ga,\kappa\>m$ in order to keep~$S$ in
$S^3\sminus K$; details are left to the reader.)

Defining $G_{L}:=G\cap S_{L}$ and letting $\bar S:=G_{L}\cup S^\mone$,
we have $\partial\bar S=\Ga$. On the other hand,
the twisted boundary
\[*
\partial^\Go\bar S=(\partial^\Go S^\mone+\Ga^\mone)+(\partial^\Go G_{L}-\Ga^\mone)
 =\tau\bigl(\dif(\Go)(\Ga)\bigl)-\tau\bigl(\Theta^\mone(\Ga)\bigr)
\]
is given by \autoref{cor.cylinder} and~\eqref{eq.dS},
and the statements follow from~\eqref{eq.A=rel}.
\endproof

\corollary[of \autoref{lem.d} and~\eqref{eq.H.twisted}]\label{cor.H.twisted}
There are canonical, up to multiplication by
integral powers of $\Go_i$'s,
isomorphisms
\[*
H_1^\Go(S_L)=H^1(F)/\Im\E(\Go),\qquad
H_1^\Go(S_L\sminus K)=H^1(F)/\E(\Go)(\Ker\kappa).
\def\qedsymbol{\donesymbol}\pushQED{\qed}\qedhere
\]
\endcorollary

\emph{Proof of} \autoref{th.C-complex}. 
If $\kappa=0$, then $K$ bounds a Seifert surface disjoint from~$F$ and,
hence, $K/L\equiv0$, which agrees with the statement of the theorem.

Therefore, till the rest of the proof we assume that $\kappa\ne0$. Assume
also that $F$ is connected, so that we can use the results of
\autoref{s.twisted} and \autoref{s.homo}.
Abbreviate $\E:=\E(\Go)$, so that $\E^*=\E(\Go\1)$ and
$\Ker\E^\perp=\Im\E^*$. Then, in view of \autoref{cor.H.twisted}, the last
two cases in the statement, as well as the finiteness of the slope in the first
case, are given by \autoref{prop.finite.slope}.
To compute this finite slope in the first case, we compare \autoref{cor.image.K} and
\autoref{lem.image.m}: if $\kappa=\E(\Ga)$, then
$\ell=-\<\alpha,\kappa\>m$ in $H_1^\Go(S_L\sminus K)$.

Finally, if $F$ is not connected, we can reduce inductively the number of
components by introducing pairs of close clasps as in \autoref{f.connected}.
\figure
\scalebox{.7}{
\begingroup%
  \makeatletter%
  \providecommand\color[2][]{%
    \errmessage{(Inkscape) Color is used for the text in Inkscape, but the package 'color.sty' is not loaded}%
    \renewcommand\color[2][]{}%
  }%
  \providecommand\transparent[1]{%
    \errmessage{(Inkscape) Transparency is used (non-zero) for the text in Inkscape, but the package 'transparent.sty' is not loaded}%
    \renewcommand\transparent[1]{}%
  }%
  \providecommand\rotatebox[2]{#2}%
  \newcommand*\fsize{\dimexpr\f@size pt\relax}%
  \newcommand*\lineheight[1]{\fontsize{\fsize}{#1\fsize}\selectfont}%
  \ifx\svgwidth\undefined%
    \setlength{\unitlength}{236.46532706bp}%
    \ifx\svgscale\undefined%
      \relax%
    \else%
      \setlength{\unitlength}{\unitlength * \real{\svgscale}}%
    \fi%
  \else%
    \setlength{\unitlength}{\svgwidth}%
  \fi%
  \global\let\svgwidth\undefined%
  \global\let\svgscale\undefined%
  \makeatother%
  \begin{picture}(1,0.41003701)%
    \lineheight{1}%
    \setlength\tabcolsep{0pt}%
    \put(0,0){\includegraphics[width=\unitlength,page=1]{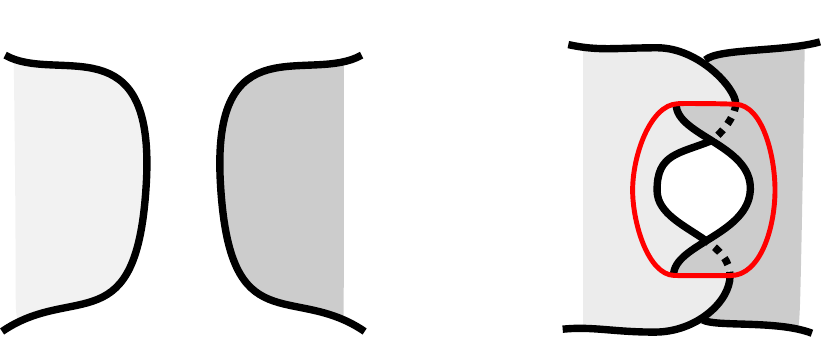}}%
    \put(0.78347656,0.29432847){\color[rgb]{1,0,0}\makebox(0,0)[lt]{\lineheight{1.25}\smash{\begin{tabular}[t]{l}$\beta$\end{tabular}}}}%
    \put(0.19761634,0.36229375){\color[rgb]{0,0,0}\makebox(0,0)[lt]{\lineheight{1.25}\smash{\begin{tabular}[t]{l}$F$\end{tabular}}}}%
    \put(0.84464582,0.38358948){\color[rgb]{0,0,0}\makebox(0,0)[lt]{\lineheight{1.25}\smash{\begin{tabular}[t]{l}$F'$\end{tabular}}}}%
  \end{picture}%
\endgroup%
}
\caption{To the left a local picture of a disconnected C-complex $F$. To the right, the complex $F'$, obtained by adding a pair of close clasps to $F$. We have $H_1(F';\Z)=H_1(F;\Z)\oplus\Z\Gb$.}
\label{f.connected}
\endfigure
If $F'$ is
obtained from~$F$ by introducing one such pair, connecting two distinct
components, then $H_1(F';\Z)=H_1(F;\Z)\oplus\Z\Gb$, where $\Gb$ is a small
proper loop running through the two clasps,
and,
extending the existing pair of dual bases by $\Gb\in H_1(F)$ and
$\Gb^*\in H^1(F)$,
the other data are
\[*
\Theta^{\prime\Ge}=\Theta^\Ge\oplus[0],\qquad
\kappa'=\kappa\oplus[0].
\]
Obviously, this modification does not affect the result of the computation.
\qed

\let\.\DOTaccent
\def\cprime{$'$}
\bibliographystyle{amsalpha}
\bibliography{bibliosplice}

\end{document}